\documentclass[oneside, reqno]{amsart}

\usepackage{xypic}
\input xy
\xyoption{all}
\usepackage{epsfig}
\usepackage{amsthm}
\usepackage{amssymb}
\usepackage{amsmath}
\usepackage{amscd}
\usepackage{color}
\usepackage[T1]{fontenc}
\usepackage{stackengine}
\stackMath

\newcommand{\bg}{\begin{equation}}
\newcommand{\ed}{\end{equation}}
\newcommand{\bga}{\begin{eqnarray}}
\newcommand{\eda}{\end{eqnarray}}
\newcommand{\pf}{\textbf{Proof:\ }}

\def\cbdu{\par{\raggedleft$\Box$\par}}

\newtheorem {Theorem}  {Theorem}

\numberwithin{Theorem}{section}

\newtheorem {Lemma}[Theorem]  {Lemma}
\newtheorem {Proposition}[Theorem]{Proposition}
\theoremstyle{definition}
\newtheorem{Definition}[Theorem]{Definition}
\theoremstyle{remark}
\newtheorem{Remark}[Theorem]{\bf Remark}

\expandafter\chardef\csname pre amssym.def
at\endcsname=\the\catcode`\@ \catcode`\@=11
\def\undefine#1{\let#1\undefined}
\def\newsymbol#1#2#3#4#5{\let\next@\relax
 \ifnum#2=\@ne\let\next@\msafam@\else
 \ifnum#2=\tw@\let\next@\msbfam@\fi\fi
 \mathchardef#1="#3\next@#4#5}
\def\mathhexbox@#1#2#3{\relax
 \ifmmode\mathpalette{}{\m@th\mathchar"#1#2#3}%
 \else\leavevmode\hbox{$\m@th\mathchar"#1#2#3$}\fi}
\def\hexnumber@#1{\ifcase#1 0\or 1\or 2\or 3\or 4\or 5\or 6\or 7\or 8\or
 9\or A\or B\or C\or D\or E\or F\fi}

\font\teneufm=eufm10 \font\seveneufm=eufm7 \font\fiveeufm=eufm5
\newfam\eufmfam
\textfont\eufmfam=\teneufm \scriptfont\eufmfam=\seveneufm
\scriptscriptfont\eufmfam=\fiveeufm

\catcode`\@=\csname pre amssym.def at\endcsname

\newcounter{remark}
\renewcommand{\div}{\mbox{div}}

\def  \12  {{\frac{1}{2}}}

\def\build#1_#2^#3{\mathrel{\mathop{\kern 0pt#1}\limits_{#2}^{#3}}}

\numberwithin{equation}{section}

\begin{document}

\title[Non-uniqueness of Hall-MHD]{Non-unique weak solutions in Leray-Hopf class for the 3D Hall-MHD system}


\author [Mimi Dai]{Mimi Dai}

\address{Department of Mathematics, Statistics, and Computer Science, University of Illinois at Chicago, Chicago, IL 60607, USA}
\email{mdai@uic.edu}

\thanks{The author was partially supported by NSF grants
DMS--1815069 and DMS--2009422.}





\begin{abstract}
Non-unique weak solutions in Leray-Hopf class are constructed for the three dimensional magneto-hydrodynamics with Hall effect. 
We adapt the widely appreciated convex integration framework developed in a recent work of Buckmaster and Vicol \cite{BV} for the Navier-Stokes equation, and with deep roots in a sequence of breakthrough papers for the Euler equation. 

\bigskip


\hspace{0.02cm}CLASSIFICATION CODE: 76D03, 76W05, 35Q35, 35D35.
\end{abstract}

\maketitle

\section{Introduction}

To capture the fast process of the magnetic reconnection phenomena in plasma physics, the following  model of the 
incompressible magneto-hydrodynamics (MHD) with Hall effect
\begin{equation}\label{HMHD}
\begin{split}
u_t+u\cdot\nabla u-B\cdot\nabla B+\nabla p=&\ \Delta u,\\
B_t+u\cdot\nabla B-B\cdot\nabla u+\zeta\nabla\times((\nabla\times B)\times B)=&\ \Delta B,\\
\nabla \cdot u=&\ 0, 
\end{split}
\end{equation}
was proposed by astrophysicists. 
In (\ref{HMHD}), $u$, $p$ and $B$ represent the fluid velocity field, the scalar pressure, and the magnetic field, respectively; they are the unknown functions on the spacial-time domain $\Omega\times [0,\infty)$. In the present paper, we take $\Omega=\mathbb T^3$. 
The parameter $\zeta$ in front of the Hall term indicates the strength of the Hall effect.
For mathematical study on this model during the last few decades, we refer to  \cite{ADFL, CL, CS, CWW, CW, Dai1, DL, DS} and references therein.

We notice that system (\ref{HMHD}) with $\zeta=0$ is the usual MHD model. In this case, one also observes that the magnetic field equation is essentially linear in $B$, while the velocity equation is obviously the Navier-Stokes equation (NSE) with a force term. Due to the linear feature of the magnetic field equation, it is expected that the properties of solutions to the MHD system do not seriously deviate from those of the solutions to the NSE. In fact, a vast amount of work for the MHD and the NSE have shown this consistence. 

However, for the Hall MHD system (\ref{HMHD}) with $\zeta>0$, the situation is drastically different, comparing to the usual MHD system. On one hand, the equation of $B$ is nonlinear with a strong nonlinear Hall term which is actually more singular than $u\cdot \nabla u$ in the NSE; on the other hand, a natural scaling does not exist for the Hall MHD system, while the MHD system shares the same natural scaling as for the NSE. More discussion on the scaling analysis will be provided at a later point. Due to the obvious difference of the two systems, a natural question is that: 
how does the presence of the Hall term change the behavior of solutions? 
Since the Hall term is more singular than other nonlinear terms in the system, one expectation is that it is probably more approachable to construct wild solutions and to discover severe ill-posedness for the Hall-MHD system. Searching wild solutions and justifying ill-posedness for fluid equations remains mathematically interesting and physically important before one can give an affirmative answer to the global regularity problem of these equations. 

As for the 3D NSE,  Leray's conjecture regarding the appearance of singularity at finite time has been a long-standing open problem; the uniqueness of Leray-Hopf weak solutions is not known either. Since the time of these problems raised in 1930s, much effort has been taken to tackle them from the negative side in the means of constructing blow-up solutions, wild solutions, or wild data-to-solution maps. Wild solutions for the Euler equation were first constructed in \cite{Sch, Sh1, Sh2}. 
In \cite{JS1}, Jia and \v{S}ver\'ak showed non-uniqueness of Leray-Hopf weak solutions in $L^\infty(L^{3,\infty})$ with the assumption that certain spectral condition holds for a linearized Navier-Stokes operator.
For an averaged NSE with modified nonlinear term, Tao constructed a smooth solution which blows up in finite time in \cite{Tao}; moreover, the author proposed a program for adapting the blowup construction to the true NSE.
Very recently, in the groundbreaking paper \cite{BV}, Buckmaster and Vicol constructed nontrivial weak solutions for the 3D NSE by developing the convex integration scheme using intermittent Beltrami flows, which leads to non-uniqueness of weak solutions.  It is a significant progress towards settling the problem of non-uniqueness of Leray-Hopf weak solutions, although the weak solutions constructed there belong to the space $C^0(0,T; H^\beta(\mathbb T^3))$ for a very small number $\beta>0$. 

The convex integration method developed in \cite{BV} dates back to \cite{MSv} and a sequence of breakthrough work for the Euler equation in the last decade. It was first introduced by De Lellis and Sz\'ekelyhidi in \cite{DLS1, DLS3, DLS2} to 
study the non-uniqueness of weak solutions and the existence of dissipative continuous solutions for the Euler equation. The framework of convex integration was further developed in \cite{BDLIS, BDLS, DSz, Is0} and eventually leads to a complete resolution of the second half of Onsager's conjecture \cite{On} by Isett \cite{Is}, and Buckmaster, De Lellis, Sz\'ekelyhidi and V. Vicol \cite{BDLSV}.

Back to the dissipative equations, as mentioned above, the non-uniqueness of Leray-Hopf weak solutions to the 3D NSE is still open. Following the convex integration method in \cite{BV}, one may expect to construct non-trivial solutions in $C^0(H^\beta)$ for $\beta<1/2$ and close enough to $1/2$; while crossing $1/2$ spacial regularity would be a major barrier. The reason is that $\dot H^{1/2}$ is critical for the 3D NSE, in which the regularity implies uniqueness. When the dissipation is weak, as for the hyperviscous Navier-Stokes equation with fractional Laplacian $(-\Delta )^\theta$ with $\theta\in(0,1/5)$ in \cite{CDLDR}, Colombo, De Lellis and De Rosa showed the non-uniqueness of Leray weak solutions, that is, solutions with finite energy and in the space $C^0(H^\theta)$. 
The result was extended to the case with $\theta<\frac13$ by De Rosa \cite{DR}.

Regarding the hyperviscous NSE with $\theta<5/4$, adapting the convex integration techniques of \cite{BV}, Luo and Titi in \cite{LT} established the non-uniqueness of weak solutions by slightly simplifying the original construction of Buckmaster and Vicol. In another work \cite{BCV}, Buckmaster, Colombo, and Vicol constructed  wild solutions for the 3D NSE, whose singular set in time has Hausdorff dimension strictly less than 1. Moreover, the result holds for the hyperviscous NSE with $\theta<5/4$ as well. Thus, along with the uniqueness result for $\theta\geq 5/4$ by Lions \cite{Lions}, the work of \cite{LT} and \cite{BCV} indicates the well-posedness criticality of the exponent $\theta=5/4$. 

Other wild solutions for the Navier-Stokes equations were also constructed in \cite{BCV, CL1, CL2, CL3}. Nonuniqueness of weak solutions was studied for other fluid equations as well, see \cite{BSV, BMS, No}. For more detailed background on recent development of convex integration method for fluid equations, the reader can consult the survey papers \cite{BV19, BV21}.

The main purpose of this paper is to address the problem of non-uniqueness of weak solutions in Leray-Hopf space for the Hall MHD system (\ref{HMHD}) with $\zeta>0$. A scaling analysis will be helpful to demonstrate why it is approachable to study this problem by adapting the convex integration techniques. We first look at the MHD system, that is (\ref{HMHD}) with $\zeta=0$.  
If the triplet $(u(x,t), B(x,t), p(x,t))$ solves the MHD system with data
$(u_0(x), B_0(x))$, the following scaled functions 
\begin{equation}\label{scal1}
u_\lambda=\lambda u(\lambda x,\lambda^2t), \ \ B_\lambda=\lambda B(\lambda x,\lambda^2t), \ \ p_\lambda=\lambda^2 p(\lambda x,\lambda^2t),
\end{equation}
solve the MHD system as well with scaled data $(\lambda u_0(\lambda x), \lambda B_0(\lambda x))$. In the case of vanishing magnetic field $B$, such scaling holds for the NSE. Under scaling (\ref{scal1}), the space $H^{1/2}(\mathbb T^3)\times H^{1/2}(\mathbb T^3)$ is critical for the 3D MHD system. It is known that regularity and hence uniqueness holds in subcritical spaces $H^{s}$ with $s>1/2$. Since the MHD system with zero magnetic field reduces to the NSE, the non-uniqueness result of the 3D NSE in \cite{BV} immediately provides a proof of non-uniqueness of weak solutions for the 3D MHD system. Similarly as for the 3D NSE, the uniqueness of Leray-Hopf solutions to the 3D MHD remains an open problem. The attempt to construct non-unique Leray-Hopf solutions via the convex integration method might not succeed since the criticality of $1/2$ spacial regularity would be a crucial obstacle to overcome. 

Now we turn to the Hall MHD system (\ref{HMHD}) with $\zeta>0$,  a natural scaling no longer holds due to the presence of the Hall term $\nabla\times((\nabla\times B)\times B)$. One can see that the Hall term is more singular than other nonlinear terms in the system and the most singular one in the magnetic equation.  This motivates us to consider the magnetic equation with vanishing velocity field as the first step. Thus we analyze the so-called Hall equation (also referred as the electron MHD),  
\begin{equation}\label{Hall}
B_t+\nabla\times((\nabla\times B)\times B)=\Delta B
\end{equation}
which has the natural scaling
\begin{equation}\label{scal2}
B_\lambda=B(\lambda x,\lambda^2t).
\end{equation}

We observe that if $\nabla\cdot B(x,0)=0$,  $\nabla\cdot B(x,t)=0$ holds for all time $t>0$.
The basic energy law for the Hall equation (\ref{Hall}) is 
\begin{equation}\label{energy-e}
\frac12\frac{d}{dt}\|B(t)\|_{L^2}^2+\|\nabla B\|_{L^2}^2=0.
\end{equation}
A Leray-Hopf type of weak solution to (\ref{Hall}) is a function $B\in L^\infty(L^2)\cap L^2(H^1)$ which satisfies (\ref{Hall}) and (\ref{energy-e}) in the distributional sense. On the other hand, under scaling (\ref{scal2}), the Sobolev space $\dot H^{3/2}$ (the same as $H^{3/2}$ on periodic domains) is critical for (\ref{Hall}) in three dimensions. One can expect global regularity of solution in $\dot H^{3/2}$ and uniqueness in the space as a consequence. Since the critical index $3/2$ of regularity is larger than the Leray-Hopf weak solution regularity index $1$, non-uniqueness of Leray-Hopf weak solutions in $C^0(H^1)$ constructed by the convex integration method would not contradict 
with anything according to the scaling properties.

Inspired by the aforementioned analysis, we adapt the convex integration scheme to the Hall equation (\ref{Hall}) and establish the first main result as follows.

\begin{Theorem}\label{thm-h}
For any nonnegative smooth function $E(t):[0,T]\to \mathbb R_{\geq 0}$, there exists a weak solution $B\in L^\infty([0,T]; L^2(\mathbb T^3))\cap C^0([0,T]; H^1(\mathbb T^3))$ to the Hall equation (\ref{Hall}), such that
\[\int_{\mathbb T^3}|\nabla\times B|^2\mathrm dx=E(t), \ \ t\in[0,T].\]
\end{Theorem}

The statement implies non-uniqueness of weak solutions to the Hall equation (\ref{Hall}) in Leray-Hopf class. Indeed, the vorticity of the weak solutions can have any nonnegative energy profiles and thus a constant (in particular, zero) is not the only weak solution.

Concerning the strategy to prove Theorem \ref{thm-h}, we take the curl of the Hall equation and apply the convex integration method to the resulted equation of the current density $J=\nabla\times B$. Section \ref{sec-hall} will be devoted to this purpose.

Once we have the convex integration scheme for the Hall equation, we turn to the coupled Hall MHD system. At each level of the convex integration which produces $B_q$, we solve the velocity field equation -- the NSE with the Lorentz force $(B_q\cdot\nabla) B_q$. We show that there exists a Leray-Hopf weak solution $u_q$ to the NSE based on the estimates on $B_q$. In the end, we illustrate that the sequence $\{u_q, B_q\}$ converges to a pair of functions $\{u, B\}$ which is a weak solution of the Hall MHD system (\ref{HMHD}). Therefore, we are able to prove the second main result stated below.

\begin{Theorem}\label{thm}
For any nonnegative smooth function $E(t):[0,T]\to \mathbb R_{\geq 0}$, there exists a weak solution  $(u, B)$ to the Hall MHD system (\ref{HMHD}) with $\zeta>0$ on $[0,T]$, such that, $u,B\in L^\infty([0,T]; L^2(\mathbb T^3))\cap L^2([0,T]; H^1(\mathbb T^3))$ and
\[\int_{\mathbb T^3}|\nabla\times B|^2\mathrm dx=E(t), \ \ t\in[0,T].\]
\end{Theorem}

Analogously, Theorem \ref{thm} suggests $(0, 0, p)$ is not the only Leray-Hopf weak solution of (\ref{HMHD}). Thus we provide a construction of non-unique Leray-Hopf weak solutions for the 3D Hall MHD system. The proof of Theorem \ref{thm} will be laid out in Section \ref{sec-hmhd}.

In the regime of applying convex integration techniques to the MHD system, the two recent articles \cite{BBV} and \cite{FLS} obtained astonishing results. In \cite{BBV}, Beekie, Buckmaster and Vicol constructed finite energy weak solutions to the ideal MHD system whose magnetic helicity is not conserved. Consequently, it shows the existence of finite energy weak solutions to the ideal MHD which cannot be obtained in the zero viscosity-resistivity limit.  It indicates that the ideal-MHD version of Taylor's conjecture is false. In \cite{FLS}, Faraco, Lindderg and Sz\'ekelyhidi showed the existence of infinitely many bounded solutions to the 3D ideal MHD. Interestingly, for such solutions, the total energy and cross helicity are not conserved; but the magnetic helicity is preserved.  Moreover, the authors showed that the 2D ideal MHD has no nontrivial compactly supported solutions with finite energy.

We conclude this section by a few well-posedness results for the Hall MHD system. In a previous paper \cite{Dai1}, the author showed that system (\ref{HMHD}) with $\zeta>0$ is locally well-posed in the Sobolev space $H^s(\mathbb R^n)\times H^s(\mathbb R^n)$ with $s>n/2$. Eventually in \cite{Dai2}, the author established the local well-posedness of the system in $H^s(\mathbb R^n)\times H^{s+1}(\mathbb R^n)$ with $s>n/2-1$, which appears to be optimal in regards to the scaling of the NSE and scaling (\ref{scal2}) for the Hall equation.

\bigskip

\section{Preliminaries}
\label{sec-pre}

\subsection{Notation}
For the sake of brevity, we first fix some notations. We denote by:
 $A\lesssim B$ an estimate of the form $A\leq C B$ with
an absolute constant $C$; $A\sim B$ an estimate of the form $C_1
B\leq A\leq C_2 B$ with absolute constants $C_1$, $C_2$.

\medskip

\subsection {The Hall equation} To analyze the effect of the Hall term, we first consider the Hall equation, which is recalled here
\begin{equation}\label{eq-H}
\begin{split}
B_t+\nabla\times((\nabla\times B)\times B)=\Delta B.
\end{split}
\end{equation}
Note that $\nabla\cdot B(t)=0$ for all $t\geq 0$ if $\nabla\cdot B(0)=0$. It is easy to verify that a smooth solution of the Hall equation satisfies the energy identity,
\begin{equation}\notag
\frac12\frac{d}{dt}\|B(t)\|_{L^2}^2+\|\nabla B(t)\|_{L^2}^2=0.
\end{equation}

\begin{Definition}\label{weak-sol}
We say that $B$ is a Leray-Hopf weak solution of (\ref{eq-H}), if for any $\varphi\in C^\infty_c([0,T]\times \mathbb T^3)$, the following integral equation 
\begin{equation}\notag
\int_0^T\int_{\mathbb T^3}B\cdot \varphi_t+(B\otimes B):\nabla\nabla\times \varphi\,\mathrm dx\,\mathrm dt=\int_0^T\int_{\mathbb T^3}\nabla B:\nabla \varphi\,\mathrm dx\,\mathrm dt
\end{equation}
is satisfied, and 
$B$ belongs to $L^\infty(0,T; L^2(\mathbb T^3))\cap L^2(0,T; H^1(\mathbb T^3))$.
\end{Definition}

\begin{Remark}
Note that the definition is valid thanks to the vector identity
\[(\nabla\times B)\times B=\nabla\cdot(B\otimes B)-\nabla \frac{|B|^2}2\]
which holds for divergence free vector field $B$. Indeed, the integral with the nonlinear term $\nabla\times((\nabla\times B)\times B)$ can be handled as, by using integration by parts
\begin{equation}\notag
\begin{split}
&\int_0^T\int_{\mathbb T^3} \nabla\times((\nabla\times B)\times B)\cdot \varphi \,\mathrm dx\,\mathrm dt\\
=& \int_0^T\int_{\mathbb T^3} \left( \nabla\cdot(B\otimes B)-\nabla \frac{|B|^2}2\right)\cdot \nabla\times\varphi \,\mathrm dx\,\mathrm dt\\
=& -\int_0^T\int_{\mathbb T^3} (B\otimes B) : \nabla \nabla\times\varphi \,\mathrm dx\,\mathrm dt+\int_0^T\int_{\mathbb T^3} \frac{|B|^2}2  \nabla \cdot (\nabla\times\varphi) \,\mathrm dx\,\mathrm dt \\
=& -\int_0^T\int_{\mathbb T^3} (B\otimes B) : \nabla \nabla\times\varphi \,\mathrm dx\,\mathrm dt,
\end{split}
\end{equation}
where the last step follows from the fact that $ \nabla \cdot (\nabla\times\varphi)=0$.
\end{Remark}

The existence of Leray-Hopf weak solutions to (\ref{eq-H}) is trivial; for instance, it can be established by the standard Galerkin's approximating method. 

Taking curl on the Hall equation leads to 
\begin{equation}\notag
(\nabla\times B)_t+\nabla\times\nabla\times((\nabla\times B)\times B)=\Delta\nabla\times B.
\end{equation}
By introducing the vorticity of the magnetic field, current density, $J=\nabla\times B$, we give two formulations of the equation. 
The first one reads as 
\begin{equation}\label{eq-curl1}
J_t+ \nabla\times\nabla\times (J\times B)=\Delta J.
\end{equation}
By applying a few vector calculus identities, see Section \ref{sec-vec}, the current density equation 
can be formulated in a more symmetric way, namely
\begin{equation}\label{eq-curl}
J_t+\nabla\cdot (B\otimes (\nabla\times J)-(\nabla\times J)\otimes B)-2\nabla\times (J\nabla B)=\Delta J
\end{equation}
where $J\nabla B=\left(J_i\partial_1 B_i, J_i\partial_2 B_i, J_i\partial_3 B_i\right)$.
In the rest of the paper, we will need to refer to both formulations to have a more complete vision of the structure of the Hall equation.

\medskip

\subsection{Leray-Hopf weak solution of the Hall-MHD}

\begin{Definition}\label{weak-sol-h}
We say that $(u, p, B)$ is a Leray-Hopf weak solution of (\ref{HMHD}), if for any $\psi,\varphi\in C^\infty_c([0,T]\times \mathbb T^3)$ with $\nabla\cdot \psi=0$, the following integral equations
\begin{equation}\notag
\int_0^T\int_{\mathbb T^3}u\cdot \psi_t+(u\otimes u):\nabla \psi-(B\otimes B):\nabla \psi\,\mathrm dx\,\mathrm dt=\int_0^T\int_{\mathbb T^3}\nabla u:\nabla \psi\,\mathrm dx\,\mathrm dt,
\end{equation}
\begin{equation}\notag
\begin{split}
\int_0^T\int_{\mathbb T^3}B\cdot \varphi_t+(u\otimes B):\nabla \varphi-(B\otimes u):\nabla \varphi+\zeta(B\otimes B):\nabla\nabla\times \varphi\,\mathrm dx\,\mathrm dt\\
=\int_0^T\int_{\mathbb T^3}\nabla B:\nabla \varphi\,\mathrm dx\,\mathrm dt
\end{split}
\end{equation}
are satisfied; 
and 
\[(u,B)\in \left(L^\infty(0,T; L^2(\mathbb T^3))\cap L^2(0,T; H^1(\mathbb T^3))\right)^2.\]
\end{Definition}
The existence of Leray-Hopf weak solutions of (\ref{HMHD}) can be found in \cite{CDL}.

\medskip

\subsection{Estimates for periodic functions and anti-derivative operator}

The following lemma regards improved H\"older's inequality for periodic functions.

\begin{Lemma}\label{le-holder} \cite{MS}
Let $\lambda\in \mathbb N$ and $f, g: \mathbb T^d\to \mathbb R$ be smooth functions. Denote $g_{\lambda}:\mathbb T^d \to \mathbb R$ by the $1/\lambda$ periodic function defined by $g_{\lambda}(x):= g(\lambda x)$. 
Then for every $p\in [1,\infty]$, we have
\begin{equation}\notag
\|f g_{\lambda}\|_{L^p}\leq   \|f\|_{L^p}\|g\|_{L^p}+ \frac{C_p}{\lambda^{1/p}} \|f\|_{C^1}\|g\|_{L^p}.
\end{equation}
\end{Lemma}
We note that 
\begin{equation}\notag
\|D^k g_{\lambda}\|_{L^p(\mathbb T^d)} =\lambda^k  \|D^k g\|_{L^p(\mathbb T^d)}, \ \ \forall \ \ k\in \mathbb N, \ \ \forall \ p\in [1,\infty],
\end{equation}
and in particular,
\begin{equation}\notag
\|g_{\lambda}\|_{L^p(\mathbb T^d)} = \|g\|_{L^p(\mathbb T^d)}, \ \ \forall \ p\in [1,\infty].
\end{equation}

A type of commutator estimate for periodic functions is introduced below.
\begin{Lemma}\label{le-comm} \cite{BV}
Assume $\kappa\geq 1$, $p\in(1,2]$ and $L\in\mathbb N$ is sufficiently large. Let function $a\in C^L(\mathbb T^3)$ be such that there exists $1\leq \lambda\leq \kappa$ and $C_a>0$ with
\[\|D^ja\|_{L^\infty}\leq C_a\lambda^j\]
for all $0\leq j\leq L$. Assume in addition that $\int_{\mathbb T^3}a(x)\mathbb P_{\geq \kappa}f(x)\, \mathrm dx=0$. Then the estimate
\[\||\nabla|^{-1}(a\mathbb P_{\geq \kappa}f)\|_{L^p}\lesssim C_a\left(1+\frac{\lambda^L}{\kappa^{L-2}}\right)\frac{\|f\|_{L^p}}{\kappa}\]
holds for any $f\in L^p(\mathbb T^3)$, where the implicit constant depends on $p$ and $L$.
\end{Lemma}
We also introduce an estimate for the symmetric anti-divergence operator.
\begin{Lemma}\label{le-anti}\cite{DLS1}
There exists a linear operator $\mathcal R$ of order $-1$, such that
\begin{equation}\notag
\nabla\cdot \mathcal R(u)=u-\displaystyle\stackinset{c}{}{c}{}{-\mkern4mu}{\displaystyle\int_{\mathbb T^3}} u\, \mathrm dx.
\end{equation}
It satisfies the Calderon-Zygmund and Schauder estimates, for $1<p<\infty$,
\begin{equation}\notag
\|\mathcal R\|_{L^p\to W^{1,p}}\lesssim 1, \ \ \|\mathcal R\|_{C^0\to C^0}\lesssim 1, \ \ 
\|\mathcal R\mathbb P_{\neq0}u\|_{L^p}\lesssim \||\nabla|^{-1}\mathbb P_{\neq0}u\|_{L^p}.
\end{equation}
\end{Lemma}

\bigskip

\section{The Hall equation and intermittent Beltrami flows}
\label{sec-hall}

In this part, we analyze the structure of the equation of the current density $J=\nabla\times B$ and lay out the intermittent Beltrami flows introduced in \cite{BV}. The analysis will reveal the fact that the equation of the current density is analogous to the NSE near the intermittent Beltrami flows.

\medskip

\subsection{Analyzing the equation}
If we apply the convex integration scheme directly to equation (\ref{eq-curl1}), we would consider the approximating equation
\begin{equation}\label{eq-curl-q}
\partial_t J_q+\nabla\times\nabla\times(J_q\times B_q)=\Delta J_q+\nabla\times\nabla\times M_q,
\end{equation}
with $J_q=\nabla\times B_q$, and $M_q$ being certain vector with the property that $M_q\to 0$ in an appropriate sense as $q\to\infty$. The main idea would be to construct building blocks for the increments $v_{q+1}=B_{q+1}-B_q$ and $w_{q+1}=J_{q+1}-J_q$, which give rise to a new pair $(B_{q+1}, J_{q+1})$ and consequently a new vector $M_{q+1}$ according to equation (\ref{eq-curl-q}) at the level of $q+1$. Most importantly, the construction should be designed in such a way that: at level $q+1$, the major contribution of nonlinear interaction to the new vector $M_{q+1}$ cancels $M_q$; and hence the sequence $\{M_q\}$ converges to zero eventually. 

However, we realize that it has certain advantages to apply the convex integration scheme to the slightly more symmetric equation (\ref{eq-curl}). In fact, we will work with the approximating form of 
(\ref{eq-curl})
\begin{equation}\label{eq-div-q}
\partial_t J_q+\nabla\cdot (B_q\otimes (\nabla\times J_q)-(\nabla\times J_q)\otimes B_q)-2\nabla \times(J_q\nabla B_q)=\Delta J_q+\nabla\cdot R_q
\end{equation}
where $R_q$ is recognized as an error stress tensor. Note that (\ref{eq-div-q}) is formally equivalent to 
\begin{equation}\label{eq-div-q-alt}
\partial_t J_q+\nabla\times\nabla\times(J_q\times B_q)=\Delta J_q+\nabla\cdot R_q.
\end{equation}
Both (\ref{eq-div-q}) and (\ref{eq-div-q-alt}) will be used in the analysis below. 
To keep the algebra as simple as possible, we may refer (\ref{eq-div-q-alt}) rather than (\ref{eq-div-q}).

 The main element is that we need to design building blocks for the increments $v_{q+1}$ and $w_{q+1}$, which in turn yield the triplet $(B_{q+1}, J_{q+1}, R_{q+1})$ with the property: the significantly large part of $R_{q+1}$ from the nonlinear interaction represented by
$B\otimes (\nabla\times J)-(\nabla\times J)\otimes B$ 
cancels the previous level of stress tensor $R_q$. A crucial observation is that:
\begin{itemize}
\item if we take $B=W(x)$ as the Beltrami wave defined in Section \ref{sec-blocks} and $J=\nabla\times W(x)=\lambda W(x)$, then we can verify 
\[ \nabla\cdot (B\otimes (\nabla\times J)-(\nabla\times J)\otimes B)-2\nabla \times (J\nabla B)=0;\]
\item if we take $B=\eta(x,t)W(x)$ as the intermittent Beltrami wave and $J=\nabla\times (\eta(x,t)W(x))$, with an approriate choice of $\eta(t,x)$ defined in Section \ref{sec-blocks}, one can make sure the difference 
\begin{equation}\label{diff}
\left[\nabla\cdot (B\otimes (\nabla\times J)-(\nabla\times J)\otimes B)-2\nabla \times(J\nabla B)-\nabla\Pi\right]-\nabla\cdot (J\otimes J)
\end{equation}
is small for some pressure term $\Pi$.
\end{itemize}
This indicates that, the stationary Beltrami wave is a solution of the Hall equation; while near certain intermittent Beltrami waves, equation (\ref{eq-curl}) is ``close'' to 
\[J_t+\nabla\cdot (J\otimes J)+\nabla \Pi=\Delta J\]
which is the NSE. 
An important motivation we obain is that an analogous construction scheme by using the convex integration method as for the NSE in \cite{BV} would possibly lead to the non-uniqueness of weak solutions of equation (\ref{eq-curl}) with $J\in C^0(0,T; H^\beta)$ for a small $\beta>0$; hence it implies $B\in C^0(0,T; \dot H^1)$ since $B$ is divergence free. Of course, in our case, two functions $J$ and $B$ are simultaneously involved in the construction; and the relation $J=\nabla\times B=\nabla\times (\eta W)$ will generate many error terms. On the other hand, it is also crucial to determine how to apply the important geometric lemma in the current context.

We will describe the convex integration scheme in detail for equation (\ref{eq-curl}) by considering its approximation sequence (\ref{eq-div-q}) in the following section.

\medskip

\subsection{Building blocks}
\label{sec-blocks}
We adapt the construction idea of \cite{BV} using intermittent Beltrami flows. While we have to keep in mind that, rather than dealing with one function satisfying the NSE, we deal with the pair $(B, J)$ with $J=\nabla\times B$ satisfying (\ref{eq-curl}) in our context.

We first fix $\xi$, $A_\xi$, $B_\xi$, and $a_\xi$ as defined in \cite{BV}:
\[\xi\in \mathbb S^2\cap \mathbb Q^3, \ \ A_\xi \in \mathbb S^2\cap \mathbb Q^3, \ \ A_\xi\cdot \xi=0, \ \ A_\xi=A_{-\xi},\]
\[a_\xi\in \mathbb C, \ \ \bar a_\xi=a_{-\xi},\]
\[B_\xi=\frac{1}{\sqrt 2}\left(A_\xi+i\xi\times A_\xi\right).
\]
The stationary Beltrami wave is taken as
\[W(x) =\sum_{\xi\in \Lambda}a_{\xi}W_\xi:=\sum_{\xi\in \Lambda}a_{\xi}B_\xi e^{i\lambda\xi\cdot x},\]
where $\Lambda$ is a given finite subset of $\mathbb S^2$ such that $\Lambda=-\Lambda$, and $\lambda$ is an integer such that $\lambda\Lambda\subset \mathbb Z^3.$
One can verify that $W(x)$ is real-valued and satifies
\[\nabla\cdot W=0,  \ \ \nabla\times W=\lambda W, \ \ \nabla\cdot (W\otimes W)=\nabla\frac{|W|^2}{2},\]
\[\displaystyle\stackinset{c}{}{c}{}{-\mkern4mu}{\displaystyle\int_{\mathbb T^3}} W\otimes W\, \mathrm dx=\frac12\sum_{\xi\in \Lambda}|a_\xi|^2(\mathrm{Id}-\xi\otimes\xi).\]

\begin{Lemma}\label{le-geo}\cite{BV}
For any $N\in\mathbb N$, we can find $\varepsilon_\gamma>0$ and $\lambda>1$ with the following property. Let $B_{\varepsilon_\gamma}(\mathrm {Id})$ be the ball of symmetric $3\times 3$ matrices, centered at $\mathrm{Id}$ of radius $\varepsilon_\gamma$. There exists pairwise disjoint subsets
\[\Lambda_\alpha\subset \mathbb S^2\cap\mathbb Q^3, \ \ \alpha\in\{1, ..., N\},\]
with $\lambda \Lambda_{\alpha}\in \mathbb Z^3$, and smooth positive functions
\[\gamma_{\xi}^{\alpha}\in C^\infty(B_{\varepsilon_\gamma}(\mathrm{Id})), \ \ \alpha\in\{1, ..., N\}, \ \ \xi\in\Lambda_\alpha,\]
with derivatives that are bounded independently of $\lambda$, such that:
\begin{enumerate}
\item $\xi\in\Lambda_{\alpha}$ implies $-\xi\in\Lambda_{\alpha}$ and $\gamma_{\xi}^\alpha=\gamma_{-\xi}^\alpha$;
\item For each $R\in B_{\varepsilon_\gamma}(\mathrm{Id})$ we have the identity
\[R=\frac12\sum_{\xi\in\Lambda_{\alpha}}\left(\gamma_{\xi}^\alpha(R)\right)^2(\mathrm{Id}-\xi\otimes\xi).\]
\end{enumerate}

\end{Lemma}

Next we describe the intermittent Beltrami flows by adding oscillations to the Beltrami waves. We start with the Dirichlet kernel $D_n$
\begin{equation}\notag
D_n(x)=\sum_{\xi=-n}^{n}e^{ix\xi}=\frac{\sin((n+\frac12)x)}{\sin(\frac x2)}
\end{equation}
which satisfies for $p>1$
\[\|D_n\|_{L^p}\sim n^{1-\frac1p}.\]
We define the lattice cube
\[\Omega_r:=\left\{\xi=(j,k,l):j,k,l\in\{-r,..., r\}\right\}\]
and the 3D normalized Dirichlet kernel
\[D_r(x):=\frac{1}{(2r+1)^{\frac32}}\sum_{\xi\in\Omega_r}e^{i\xi\cdot x}\]
satisfying
\begin{equation}\notag
\|D_r\|_{L^2}^2=(2\pi)^3, \ \ \|D_r\|_{L^p}\lesssim r^{\frac32-\frac3p}, \ \ p>1,
\end{equation}
where the implicit constant depends only on $p$. The parameter $r$ refers to the number of frequencies along edges of $\Omega_r$.

We shall define a directed and rescaled periodic Dirichlet kernel with period $\left(\mathbb T/(\lambda\sigma)\right)^3$. The small constant $\sigma$ is chosen such that $\lambda \sigma \in \mathbb N$ which parameterizes the spacing between frequencies; and $\sigma r\ll 1$. We fix an integer $N_0\geq 1$ such that 
\[\{N_0\xi, N_0A_\xi, N_0\xi\times A_\xi\}\subset N_0\mathbb S^2\cap\mathbb Z^3\]
for all $\xi\in \Lambda_\alpha$ and $\alpha\in \{1, ..., N\}.$ We also introduce a parameter $\mu\in (\lambda, \lambda^2)$, which adjusts the temporal oscillation. It is then ready to define the modified Dirichlet kernel
\begin{equation}\label{eta}
\eta_{\xi, \lambda,\sigma, r, \mu}(x,t)=D_r\left(\lambda\sigma N_0(\xi\cdot x+\mu t), \lambda\sigma N_0 A_\xi\cdot x, \lambda\sigma N_0(\xi\times A_\xi)\cdot x\right)
\end{equation}
for $\xi\in \Lambda_\alpha^+$; while $\eta_{\xi, \lambda,\sigma, r, \mu}(x,t)=\eta_{-\xi, \lambda,\sigma, r, \mu}(x,t)$ for $\xi\in \Lambda_\alpha^-$. We take the short notation $\eta_{\xi}(x,t)=\eta_{\xi,\lambda,\sigma,r,\mu}(x,t)$. It is important to notice that
\begin{equation}\label{eta-1}
\mu^{-1}\partial_t\eta_{\xi}(x,t)=\pm(\xi\cdot\nabla)\eta_{\xi}(x,t), \ \forall \xi\in \Lambda_{\alpha}^{\pm},
\end{equation}
which is the crucial identity used to design temporal oscillation in the increments later.

One also observe that 
\begin{equation}\label{eta-norm}
\displaystyle\stackinset{c}{}{c}{}{-\mkern4mu}{\displaystyle\int_{\mathbb T^3}} \eta_{\xi}^2(x,t)\, \mathrm dx = \displaystyle\stackinset{c}{}{c}{}{-\mkern4mu}{\displaystyle\int_{\mathbb T^3}} D_r^2(x)\, \mathrm dx=1, \ \ 
\|\eta_{\xi}(\cdot, t)\|_{L^p}=\|D_{r}\|_{L^p}\lesssim r^{\frac32-\frac3p},
\end{equation}
for all $1<p\leq \infty$.

Now we are ready to introduce the intermittent wave $\mathbb W_{\xi}$:
\begin{equation}\label{int-wave}
\mathbb W_{\xi}(x,t)=\eta_{\xi}(x,t)B_\xi e^{i\lambda\xi\cdot x}.
\end{equation}
It is worth to point out that $\mathbb W_{\xi}$ is supported on certain frequencies. Indeed, we have
\begin{equation}\notag
\begin{split}
\mathbb P_{\leq 2\lambda\sigma rN_0}\eta_{\xi}=&\ \eta_{\xi},\\
\mathbb P_{\leq 2\lambda}\mathbb P_{\geq \lambda/2} \mathbb W_{\xi}=&\ \mathbb W_{\xi},\\
\mathbb P_{\leq 4\lambda}\mathbb P_{\geq c_0\lambda} \left( \mathbb W_{\xi}\otimes  \mathbb W_{\xi'}\right)=& \ \mathbb W_{\xi}\otimes  \mathbb W_{\xi'}
\end{split}
\end{equation}
where $c_0$ is a small constant and $\xi'\neq -\xi$. 

Another important fact regarding $\mathbb W_{\xi}$ is given by
\begin{equation}\notag
\nabla\cdot(\mathbb W_{\xi}\otimes \mathbb W_{-\xi}+\mathbb W_{-\xi}\otimes \mathbb W_{\xi})
=\nabla \eta_{\xi}^2-(\xi\cdot\nabla)\eta_{\xi}^2\xi=\nabla \eta_{\xi}^2-\frac{\xi}{\mu}\partial_t\eta_{\xi}^2.
\end{equation}
It is the main motivation that we need to include the temporal oscillation $w_{q+1}^t$ into the construction later.

Different from the Beltrami wave $W_{\xi}(x)=B_\xi e^{i\lambda\xi\cdot x}$, the intermittent Beltrami wave $\mathbb W_{\xi}$ is not divergence free or an eigenfunction of curl, i.e.
\[\nabla\cdot \mathbb W_{\xi}\neq 0, \ \ \nabla\times \mathbb W_{\xi}\neq \lambda \mathbb W_{\xi}.\]
Instead, we have
\begin{equation}\notag
\begin{split}
\nabla\cdot \mathbb W_{\xi}=&\ \nabla\eta_{\xi}\cdot W_{\xi},\\
\nabla\times\mathbb W_{\xi}=&\ \lambda \mathbb W_{\xi}+\nabla\eta_{\xi}\times W_{\xi}.
\end{split}
\end{equation}
Parameters $\lambda, \sigma, r$, and $\mu$ will be chosen in an appropriate way such that $\nabla\eta_{\xi}\cdot W_{\xi}$ and $\nabla\eta_{\xi}\times W_{\xi}$ are sufficiently small. 

For such intermittent Beltrami waves $\mathbb W_{\xi}$ and $\Lambda_\alpha, \varepsilon_\gamma, \gamma_\xi$ as in Lemma \ref{le-geo}, we have the following geometric lemma, which is a key ingredient in the construction.

\begin{Lemma}\label{le-geo-2}\cite{BV}
Assume $a_\xi\in\mathbb C$ are constants satisfying $\bar a_\xi=a_{-\xi}$.  The vector field
\[\sum_{\alpha}\sum_{\xi\in\Lambda_\alpha}a_\xi\mathbb W_{\xi}(x)\]
is real valued. Moreover, for each matrix $R\in B_{\varepsilon_\gamma}(\mathrm{Id})$ we have
\begin{equation}\label{geo-id}
\sum_{\xi\in \Lambda_\alpha}\left(\gamma_\xi(R)\right)^2\displaystyle\stackinset{c}{}{c}{}{-\mkern4mu}{\displaystyle\int_{\mathbb T^3}} \mathbb W_{\xi}\otimes \mathbb W_{-\xi}\, \mathrm dx=
\sum_{\xi\in \Lambda_\alpha}\left(\gamma_\xi(R)\right)^2 B_{\xi}\otimes B_{-\xi}=R.
\end{equation}
\end{Lemma}

\medskip

\subsection{Analogy of equation (\ref{eq-curl}) with the NSE near intermittent Beltrami flows}
In this part, we further analyze the structure of the nonlinearity of equation (\ref{eq-curl}) by comparing it with the NSE near the intermittent Beltrami flows introduced above.
We can take the magnetic field $B$ as
\[\mathbb W_{\xi}^B=\frac{1}{\lambda}\mathbb W_\xi=\frac{1}{\lambda}\eta_{\xi} W_{\xi}.\]
An important observation is that
\[\nabla\times \mathbb W_{\xi}^B=\eta_{\xi} W_{\xi}+\frac{1}{\lambda}\nabla \eta_{\xi}\times W_{\xi}=\mathbb W_{\xi}+\frac{1}{\lambda}\nabla \eta_{\xi}\times W_{\xi}\]
and 
\[\|\frac{1}{\lambda}\nabla \eta_{\xi}\times W_{\xi}\|_{L^2}\lesssim \sigma r,\]
where the upper bound $\sigma r$ can be sufficiently small by choosing the parameters appropriately. We denote $\mathbb W_\varepsilon=:\frac{1}{\lambda}\nabla \eta_{\xi}\times W_{\xi}$ to be the small error term between $\nabla\times \mathbb W_{\xi}^B$ and $\mathbb W_{\xi}$. Thus we can naturally adapt $J=\nabla\times \mathbb W_{\xi}^B$. 

Now we show that the difference (\ref{diff}) is actually small near the intermittent Beltrami flows. 
Namely, by taking $B=\mathbb W_{\xi}^B=\frac{1}{\lambda}\mathbb W_{\xi}$ and $J=\nabla\times \mathbb W_{\xi}^B=\mathbb W_{\xi}+\mathbb W_\varepsilon$, a straight forward computation shows the difference 
\begin{equation}\notag
\begin{split}
&\left[\nabla\cdot (\mathbb W_{\xi}^B\otimes (\nabla\times J)-(\nabla\times J)\otimes \mathbb W_{\xi}^B-2\nabla\times(J\nabla \mathbb W_{\xi}^B))+\nabla \frac{|J|^2}{2}\right]-\nabla\cdot (J\otimes J)\\
=&\ \nabla\cdot \left(\frac{1}{\lambda}\mathbb W_{\xi}\otimes (\nabla\times (\mathbb W_{\xi}+\mathbb W_\varepsilon)-(\nabla\times (\mathbb W_{\xi}+\mathbb W_\varepsilon))\otimes \frac{1}{\lambda}\mathbb W_{\xi}\right)\\
&-2\nabla\times \left((\mathbb W_{\xi}+\mathbb W_\varepsilon)\nabla \mathbb W_{\xi})\right)
+\nabla \frac{|\mathbb W_{\xi}+\mathbb W_\varepsilon|^2}{2}-\nabla\cdot ((\mathbb W_{\xi}+\mathbb W_\varepsilon)\otimes (\mathbb W_{\xi}+\mathbb W_\varepsilon))\\
\sim &\ \nabla\cdot \left(\frac{1}{\lambda}\mathbb W_{\xi}\otimes (\nabla\times \mathbb W_{\xi})-(\nabla\times \mathbb W_{\xi})\otimes \frac{1}{\lambda}\mathbb W_{\xi}\right)-2\nabla\times (\mathbb W_{\xi}\nabla\mathbb W_{\xi})\\
&+\nabla \frac{|\mathbb W_{\xi}|^2}{2}-\nabla\cdot (\mathbb W_{\xi}\otimes \mathbb W_{\xi})\\
\sim &\ 0.
\end{split}
\end{equation} 
Thus, near the intermittent Beltrami flows $(B, J)=(\mathbb W_{\xi}^B, \nabla\times \mathbb W_{\xi}^B)$,  equation (\ref{eq-curl}) (the curl of the Hall equation) is indeed ``close'' to the NSE. Also, an obvious fact is that $J=\nabla\times B$ scales as the velocity field in the NSE. This is the main motivation to investigate the problem of non-uniqueness of weak solutions for the Hall-MHD system in Leray-Hopf space by adhering to what has been done for the NSE in \cite{BV}. Of course, new difficulties arise in the construction. In particular, rather than one function, involved here are a pair of functions $B$ and $J$, which are related through $J=\nabla\times B$. On the other hand, to apply the rigid geometric lemma, one has to regroup the nonlinear interactions in a suitable way such that error terms can be controlled. It is also non-trivial to determine how to introduce the temporal oscillation. In the end, to show non-uniqueness of Leray-Hopf weak solutions for the Hall-MHD system,  we need to design a scheme of combining the convex integration method for the magnetic field equation and the classical regularity theory for the NSE.
We will address all of these problems in the rest of the article.

\bigskip

\section{Convex integration for the Hall equation}
\label{sec-hall}

In this part, we adapt the convex integration method to construct Leray-Hopf weak solutions of the Hall equation with nonnegative energy profiles for the current density field. The main strategy is to design an iteration scheme for the approximating equation (\ref{eq-div-q}) illustrated in Proposition \ref{le-iterative}.

We start with fixing several parameters:
for large enough constants $a\gg 1$ and $b\gg 1$, and small enough positive constant $\beta\ll 1$, we define:
\begin{equation}\label{pp1}
\lambda_q=a^{b^q}, \ \ \delta_q=\lambda_1^{3\beta}\lambda_q^{-2\beta},
\end{equation}
\begin{equation}\label{pp2}
r=\lambda_{q+1}^{\frac34}, \ \sigma=\lambda_{q+1}^{-\frac{15}{16}}, \ \mu=\lambda_{q+1}^{\frac54}, \ \ell=\lambda_q^{-20}.
\end{equation}
It is easy to see that $\lambda_{q+1}=\lambda_q^b$.

\begin{Proposition}\label{le-iterative}
There exists an absolute constant $C>0$ and a sufficiently small parameter $\varepsilon_R$ depending on $b$ and $\beta$ such that the following inductive statement holds.
Let $(B_q, J_q, R_q)$ be a solution of the approximating equation (\ref{eq-div-q}) on $\mathbb T^3\times [0,T]$ satisfying:
\begin{equation}\label{bq-c1}
\|B_q\|_{C^1_{x,t}}\leq \lambda_q^3,
\end{equation}
\begin{equation} \label{uq-c1}
\|J_q\|_{C^1_{x,t}}\leq \lambda_q^4, 
\end{equation}
\begin{equation}\label{energy-q}
0\leq E(t)-\int_{\mathbb T^3}|J_q|^2\, \mathrm dx\leq \delta_{q+1},
\end{equation}
and 
\begin{equation}\label{energy-q1}
E(t)-\int_{\mathbb T^3}|J_q|^2\, \mathrm dx\leq \frac{\delta_{q+1}}{100} \ \ \mbox {implies} \ \ J_{q}(\cdot, t)\equiv 0 \ \mbox{and} \ R_q(\cdot, t)\equiv 0.
\end{equation}
In addition, we assume 
\begin{equation}\label{RM}
\nabla\cdot R_q=\nabla\cdot \tilde R_q+\nabla\times\nabla\times \tilde M_q+\nabla\cdot\nabla \tilde Q_q+\nabla\tilde p_{q+1}
\end{equation}
for a traceless symmetric tensor $\tilde R_q$, vector field $\tilde M_q$ and $\tilde Q_q$, and a scalar pressure function $\tilde p_{q+1}$, which satisfy
\begin{equation}\label{rq-l1}
\|\tilde R_q\|_{L^\infty(L^1)}+\|\tilde M_q\|_{L^\infty(L^1)}+\|\tilde Q_q\|_{L^\infty(L^1)}\leq \lambda_q^{-\varepsilon_R}\delta_{q+1},
\end{equation}
\begin{equation}\label{lq-c1}
\|R_q\|_{C^1_{x,t}}\leq \lambda_q^{12}.
\end{equation}
Then we can find another solution $(B_{q+1}, J_{q+1}, R_{q+1})$ of (\ref{eq-div-q}) satisfying (\ref{bq-c1})-(\ref{lq-c1}) with $q$ replaced by $q+1$. Moreover, the increments $v_{q+1}=B_{q+1}-B_q$ and $w_{q+1}=J_{q+1}-J_q$ satisfy 
\begin{equation}\label{est-increm}
\|v_{q+1}\|_{L^2}\leq C\lambda_{q+1}^{-1}\delta_{q+1}^{1/2}, \ \ \|w_{q+1}\|_{L^2}\leq C\delta_{q+1}^{1/2}.
\end{equation}
\end{Proposition}
This proposition leads to a proof of Theorem \ref{thm-h} immediately. 

\textbf {Proof of Theorem \ref{thm-h}:}
At the first step, we take $(B_0, J_0, R_0)=(0,0,0)$ which satisfies (\ref{bq-c1})-(\ref{rq-l1}), and (\ref{energy-q})-(\ref{energy-q1}) for large enough $a>0$.  For $q\geq 1$, we apply Proposition \ref{le-iterative} to obtain a sequence of approximating solutions $\{(B_q, J_q, R_q)\}$ satisfying (\ref{bq-c1})-(\ref{energy-q1}). It follows from
(\ref{est-increm}) that  
\begin{equation}\notag
\sum_{q\geq0}\|J_{q+1}-J_q\|_{L^2}= \sum_{q\geq0}\|w_{q+1}\|_{L^2}\lesssim \sum_{q\geq0}\delta_{q+1}^{1/2}<\infty.
\end{equation}
which implies the strong convergence of $J_q=\nabla\times B_q$ to a function $J$ in $C^0(0,T;L^2)$, and the strong convergence of $B_q$ to a function $B$ in $C^0(0,T;H^1)$ with $J=\nabla\times B$ and $\nabla\cdot B=0$.

While $\|\tilde R_q\|_{L^\infty(0,T; L^1)}\to 0$ and $\|\tilde M_q\|_{L^\infty(0,T; L^1)}\to 0$ as $q\to \infty$, we conclude $J$ is a weak solution of  (\ref{eq-curl}), and $B$ is a weak solution of (\ref{eq-H}); moreover, it is obvious that $B\in L^\infty(0,T; L^2(\mathbb T^3))\cap L^2(0,T; H^1(\mathbb T^3))$, since $B$ is divergence free.

\cbdu

The proof of Proposition \ref{le-iterative} will be carried out in Sections \ref{sec-design} -- \ref{sec-energy} below. 

\medskip

\subsection{Construction of the perturbation $(v_{q+1}, w_{q+1})$}
\label{sec-design}
Based on the building blocks introduced in Section \ref{sec-blocks}, we proceed to construct the perturbation $v_{q+1}=B_{q+1}-B_q$,
\[v_{q+1} :=v_{q+1}^{p}+v_{q+1}^{c}+v_{q+1}^{t}\]
where $v_{q+1}^{p}$ and $v_{q+1}^{c}$  are defined as
\begin{equation}\notag
\begin{split}
v_{q+1}^{p}=&\sum_{\xi\in \Lambda_\alpha} a_\xi\mathbb W_{\xi}^B=\lambda_{q+1}^{-1}\sum_{\xi\in \Lambda_\alpha}a_{\xi}\eta_{\xi} W_{\xi},\\
v_{q+1}^{c}=& \lambda_{q+1}^{-2}\sum_{\xi\in \Lambda_\alpha}\nabla(a_{\xi}\eta_{\xi})\times W_{\xi},\\
\end{split}
\end{equation}
while $v_{q+1}^{t}$ will be defined through $w_{q+1}^t$ later. One can verify that
\[\nabla\cdot (v_{q+1}^p+v_{q+1}^c)=\lambda_{q+1}^{-2}\sum_{\xi\in \Lambda_\alpha}\nabla\cdot(\nabla\times (a_{\xi}\eta_{\xi}W_{\xi}))=0.\]
We now define the perturbation $w_{q+1}=J_{q+1}-J_q$ as
\[w_{q+1}=w_{q+1}^p+w_{q+1}^c+w_{q+1}^t\] with 
\begin{equation}\notag
\begin{split}
w_{q+1}^{p}=&\nabla\times v_{q+1}^{p}, \ w_{q+1}^{c}=\nabla\times v_{q+1}^{c}, \\
w_{q+1}^{t}=&\mu^{-1}\sum_{\xi}\mathbb P_H\mathbb P_{\neq0}(a_\xi^2\eta^2\xi).
\end{split}
\end{equation}
In the end, we define $v_{q+1}^t$ through $w_{q+1}^t=\nabla\times v^t_{q+1}$ up to a gradient which we can take as zero. Indeed, for $v_{q+1}^t\in L^2$, we can decompose $v_{q+1}^t$ as 
\[v_{q+1}^t=v_{q+1,0}^t+\nabla \phi, \ \ \mbox {with} \ \ \nabla\cdot v_{q+1,0}^t=0.\]
In our case, we simply take $v_{q+1,0}^t$ to be $v_{q+1}^t$, since $\nabla\times \nabla\phi=0$. Thus, $\nabla\cdot v_{q+1}^t=0$ holds. Along with the fact $\nabla\cdot (v_{q+1}^p+v_{q+1}^c)=0$, we have 
\[\nabla\cdot v_{q+1}=0.\]
On the other hand, it is obvious that 
\[\nabla\cdot w_{q+1}^p=\nabla\cdot w_{q+1}^c=\nabla\cdot w_{q+1}^t=\nabla\cdot w_{q+1}=0\]
and 
\[w_{q+1}=\nabla\times v_{q+1}.\]

\medskip

\subsection{Estimates of building blocks}
\label{sec-est-blocks}

The main purpose of adding the oscillation $\eta_{\xi}$ to the Beltrami waves is to make sure the $L^1$ norm of the waves is significantly smaller than the $L^2$ norm. This can be seen in the following lemma.

\begin{Lemma}\label{le-w} \cite{BV}
The bounds 
\begin{equation}
\|\nabla^N\partial_t^K\mathbb W_{\xi}\|_{L^p}\lesssim \lambda^N(\lambda\sigma r\mu)^Kr^{\frac32-\frac3p}, \label{w-lp}\\
\end{equation}
\begin{equation}
\|\nabla^N\partial_t^K\mathbb \eta_{\xi}\|_{L^p}\lesssim (\lambda\sigma r)^N(\lambda\sigma r\mu)^Kr^{\frac32-\frac3p}  \label{eta-lp}
\end{equation}
hold for all $1<p\leq \infty$.
\end{Lemma}

We point out that, following \cite{BV}, in order to avoid a loss of derivative,  the pair $(v_q, w_q)$ at each level needs to be regularized by using standard Frieddrichs mollifiers. Moreover, the corresponding stress tensor $R_q$ is not spatially homogenous. To fix it, cutoff functions that form a partition of unity can be introduced to decompose $R_q$ into slices. The two steps involve delicate computations, which will be omitted in our presentation. Rather, we do adapt the regularization parameter  $\ell$ from the first step. We also adapt the partition of unity: let $0\leq \tilde \chi_0, \tilde \chi\leq 1$ be smooth functions supported on $[0,4]$ and $[\frac14,4]$ respectively; and $\tilde \chi_i(z)=\tilde\chi(4^{-i}z)$ satisfying
\begin{equation}\notag
\tilde\chi_0^2(z)+\sum_{i\geq 1}\tilde \chi_i(z)\equiv 1, \ \ \forall z>0.
\end{equation}
Then we define the amplitude function $a_\xi$ for the intermittent Beltrami flows as,
\begin{equation}\label{amp}
a_{\xi,i,q+1}=\rho_i^{\frac12}\chi_{i,q+1}\gamma(\xi)\left(\mathrm{Id}-\frac{R_q}{\rho_i}\right)
\end{equation}
where $\rho_i$ and $\chi_{i,q+1}$ are defined as
\begin{equation}\notag
\begin{split}
\rho_i=& \ \lambda_q^{-\varepsilon_R}\delta_{q+1}4^{i+c_0}, \ \ i\geq 1,\\
\chi_{i,q+1}(x,t)= &\ \tilde\chi_i\left(\left<\frac{R_q(x,t)}{100\lambda_q^{-\varepsilon_R}\delta_{q+1}}\right>\right).
\end{split}
\end{equation}
Here we use the notation $\left<A\right>=(1+|A|^2)^{\frac12}$ with $|\cdot|$ being the Euclidean norm of a matrix. Referring to \cite{BV}, we have 
\begin{equation}\label{est-il}
4^{\max\{i\}}\lesssim \ell^{-1}.
\end{equation}
To make sure the inequality (\ref{energy-q}) holds, we need to choose $\rho_0$ as follows,
\begin{equation}\notag
\begin{split}
\rho(t)=&\ \frac1{3|\mathbb T^3|}\left(\int_{\mathbb T^3}\chi_0^2\,\mathrm dx\right)^{-1}\max\left(E(t)-\int_{\mathbb T^3}|J_q|^2\,\mathrm dx-3\sum_{i\geq 1}\rho_i\int_{\mathbb T^3}\chi_i^2\, \mathrm dx-\frac{\delta_{q+1}}{2}, 0\right)\\
\rho_0=&\ \left((\rho^{1/2})*\varphi_{\ell}\right)^2,
\end{split}
\end{equation}
where $\varphi_{\ell}$ is the standard Friedrichs mollifer at time scale $\ell$. It was shown in \cite{BV}, such defined $\rho_0$ satisfies 
\begin{equation}\notag
\|\rho_0\|_{C^0_t}\leq 2\delta_{q+1}, \ \ \|\rho_0^{1/2}\|_{C^N_t}\lesssim \delta_{q+1}^{1/2}\ell^{-N}, \ \mbox{for} \ \ N\geq 1.
\end{equation}

Below is a collection of estimates satisfied by the amplitude function $a_\xi$.

\begin{Lemma}\label{le-a}
The following bounds hold
\begin{align}
\|\chi_i\|_{L^1}\lesssim&\ 4^{-i}, \label{est-chi}\\
\|a_{\xi}\|_{L^2}\lesssim&\ \delta_{q+1}^{\frac12}, \label{a-l2}\\
\|a_{\xi}\|_{L^\infty}\lesssim&\ \delta_{q+1}^{\frac12}\ell^{-\frac12},  \label{a-inf}\\
\|a_{\xi}\|_{L^p}\lesssim&\ \delta_{q+1}^{\frac12}\ell^{-\frac12(1-\frac1p)}, \ \ \mbox {for} \ \ p\geq 1,  \label{a-lp}\\
\|a_{\xi}\|_{C^N_{x,t}}\lesssim&\ \ell^{-N}, \ \ \mbox {for} \ \ N\geq 1. \label{a-cn}
\end{align}
\end{Lemma}
\pf
We only need to show (\ref{a-lp}), since other ones were shown in \cite{BV}.
In view of (\ref{est-chi}), we deduce
\begin{equation}\notag
\|a_\xi\|_{L^1}\lesssim  \rho_i^{\frac12} \|\chi_i\|_{L^1}\lesssim \lambda_q^{-\varepsilon_R/2}\delta_{q+1}^{\frac12}.
\end{equation}
Thus, by interpolation we obtain
\begin{equation}\notag
\|a_\xi\|_{L^p}\lesssim \|a_{\xi}\|_{L^\infty}^{\frac{p-1}p}\|a_{\xi}\|_{L^1}^{\frac1p}\lesssim \delta_{q+1}^{\frac12(1-\frac1p)}\ell^{-\frac12(1-\frac1p)}\lambda_q^{-\varepsilon_R\frac1{2p}}\delta_{q+1}^{\frac1{2p}}\lesssim \delta_{q+1}^{\frac12}\ell^{-\frac12(1-\frac1p)}.
\end{equation}

\cbdu

\medskip

\subsection{Estimates of the perturbation}

\begin{Lemma}\label{le-vq}
The increment $v_{q+1}=B_{q+1}-B_q$ satisfies the following estimates
\begin{align}
\|v_{q+1}^p\|_{L^2}\lesssim&\ \lambda_{q+1}^{-1}\delta_{q+1}^{\frac12}, \label{v-p2}\\
\|v_{q+1}^c\|_{L^2}\lesssim&\ \ell^{-1}\mu^{-1}\lambda_{q+1}^{-1}\delta_{q+1}^{\frac12}r^{\frac32}, \label{v-c2}\\
\|v_{q+1}^t\|_{L^2}\lesssim&\ \ell^{-1}\mu^{-1}(\lambda_{q+1}\sigma)^{-1}\delta_{q+1}r^{\frac32}, \label{v-t2}\\
\|v_{q+1}^p\|_{L^p}\lesssim&\ \lambda_{q+1}^{-1}\delta_{q+1}^{\frac12}\ell^{-\frac12(1-\frac1p)}r^{\frac32-\frac3p}, \ \ p\geq1, \label{v-pp}\\
\|v_{q+1}^c\|_{L^p}\lesssim&\ \lambda_{q+1}^{-1}\delta_{q+1}^{\frac12}\ell^{-\frac12(1-\frac1p)}\sigma r^{\frac52-\frac3p}, \ \ p\geq 1, \label{v-cp}\\
\|v_{q+1}^t\|_{L^p}\lesssim&\ \ell^{-1}\mu^{-1}(\lambda_{q+1}\sigma)^{-1}\delta_{q+1}r^{3-\frac3p}, \ \ p\geq 1, \label{v-tp}\\
\|v_{q+1}^p\|_{W^{1,p}}+\|v_{q+1}^c\|_{W^{1,p}}\lesssim&\ \ell^{-2}r^{\frac32-\frac3p}, \ \ p\geq 1, \label{v-w1p}\\
\|v_{q+1}^t\|_{W^{1,p}}\lesssim&\ \mu^{-1}\delta_{q+1}\ell^{-1}r^{4-\frac3p}, \ \ p\geq 1, \label{vt-w1p}\\
\|v_{q+1}^p\|_{C^N_{x,t}}+\|v_{q+1}^c\|_{C^N_{x,t}}\lesssim&\ \lambda_{q+1}^{\frac{1+5N}{2}}, \label{v-cn}\\
\|v_{q+1}^t\|_{C^N_{x,t}}\lesssim&\ \lambda_{q+1}^{\frac{3+5N}{2}}, \label{v-tn}\\
\|B_{q+1}\|_{C^N_{x,t}}\lesssim&\ \lambda_{q+1}^{\frac{3+5N}{2}}. \label{b-cn}
\end{align}
\end{Lemma}
\pf 
Adhering to Lemma \ref{le-holder}, (\ref{w-lp}), (\ref{a-lp}) and (\ref{a-cn}), we obtain for $1\leq p\leq \infty$
\begin{equation}\notag
\begin{split}
\|v_{q+1}^p\|_{L^p}\leq&\ \lambda_{q+1}^{-1}\sum_{\xi\in\Lambda_\alpha}\left\|a_{\xi}\mathbb W_{\xi}\right\|_{L^p}\\
\lesssim&\ \lambda_{q+1}^{-1}\sum_{\xi\in\Lambda_{\alpha}}\left(\|a_\xi\|_{L^p}\|\mathbb W_\xi\|_{L^p}+\lambda_{q+1}^{-\frac{1}{p}} \|a_\xi\|_{C^1_x}\|\mathbb W_\xi\|_{L^p}\right)\\
\lesssim&\ \lambda_{q+1}^{-1}\delta_{q+1}^{\frac12}\ell^{-\frac12(1-\frac1p)}r^{\frac32-\frac3p}.
\end{split}
\end{equation}
In view of Lemma \ref{le-holder}, (\ref{a-lp}), (\ref{a-cn}), and (\ref{eta-lp}), and the choice of parameters (\ref{pp1})-(\ref{pp2}), we obtain
\begin{equation}\notag
\begin{split}
\|v_{q+1}^c\|_{L^p}\leq&\ \lambda_{q+1}^{-2}\sum_{\xi\in\Lambda_\alpha}\|\nabla(a_\xi\eta_\xi)\|_{L^p}\\
\lesssim&\ \lambda_{q+1}^{-2}\sum_{\xi\in\Lambda_\alpha}\|a_\xi\nabla\eta_\xi\|_{L^p}+\lambda_{q+1}^{-2}\sum_{\xi\in\Lambda_\alpha}\|\nabla a_\xi\eta_\xi\|_{L^p}\\
\lesssim&\ \lambda_{q+1}^{-2}\sum_{\xi\in\Lambda_\alpha}\left(\|a_\xi\|_{L^p}\|\nabla\eta_\xi\|_{L^p}+\lambda_{q+1}^{-\frac1p}\|a_\xi\|_{C^1_x}\|\nabla\eta_\xi\|_{L^p}\right)\\
&+\lambda_{q+1}^{-2}\sum_{\xi\in\Lambda_\alpha}\|\nabla a_\xi\|_{L^\infty}\|\eta_\xi\|_{L^p}\\
\lesssim&\ \lambda_{q+1}^{-2}\delta_{q+1}^{\frac12}\ell^{-\frac12(1-\frac1p)}\lambda_{q+1}\sigma r^{\frac52-\frac3p}+\lambda_{q+1}^{-2}\ell^{-1}r^{\frac32-\frac3p}\\
\lesssim&\ \lambda_{q+1}^{-1}\delta_{q+1}^{\frac12}\ell^{-\frac12(1-\frac1p)}\sigma r^{\frac52-\frac3p}.
\end{split}
\end{equation}

The proof of estimates on other norms of $v_{q+1}^p$ and $v_{q+1}^c$ can be found in \cite{BV} (by multiplying each estimate the factor $\lambda_{q+1}^{-1}$). We only need to show the estimates for $v_{q+1}^t$. Recall that $v_{q+1}^t$ satisfies
\[\nabla\times v_{q+1}^t=w_{q+1}^t=\mu^{-1}\sum_{\xi\in \Lambda}\mathbb P_H\mathbb P_{\neq 0}(a_{\xi}^2\eta_{\xi}^2\xi)\]
and $v_{q+1}^t$ is divergence free. Thus by Lemma \ref{le-comm} we deduce
\begin{equation}\notag
\begin{split}
\|v_{q+1}^t\|_{L^2}\leq &\mu^{-1}\left\|\sum_{\xi\in \Lambda}\mathrm{curl}^{-1}\mathbb P_H\mathbb P_{\neq 0}(a_{\xi}^2\eta_{\xi}^2\xi)\right\|_{L^2}\\
\lesssim &\mu^{-1}\left\|\sum_{\xi\in \Lambda}\mathrm{curl}^{-1}\left(a_{\xi}^2\mathbb P_{\geq \lambda_{q+1}\sigma/2}(\eta_{\xi}^2\xi)\right)\right\|_{L^2}\\
\lesssim &\mu^{-1}\sum_{\xi\in \Lambda}\|a_{\xi}^2\|_{L^\infty}(\lambda_{q+1}\sigma)^{-1}\left(1+\frac{1}{\ell^{L}(\lambda_{q+1}\sigma)^{L-2}}\right)\|\eta_{\xi}^2\|_{L^2}\\
\lesssim &\mu^{-1}\delta_{q+1}\ell^{-1}(\lambda_{q+1}\sigma)^{-1}r^{\frac32}.
\end{split}
\end{equation}
In an analogous way, we can obtain
\begin{equation}\notag
\begin{split}
\|v_{q+1}^t\|_{L^p}\lesssim &\mu^{-1}\delta_{q+1}\ell^{-1}(\lambda_{q+1}\sigma)^{-1}r^{3-\frac3p},\\
\|v_{q+1}^t\|_{W^{1,p}}\lesssim &\mu^{-1}\delta_{q+1}\ell^{-1}r^{4-\frac3p}.
\end{split}
\end{equation}
Proof of inequality (\ref{v-tn}) can be referred to \cite{BV}; inequality (\ref{b-cn}) follows from (\ref{v-cn}) and (\ref{v-tn}).

\cbdu

\begin{Lemma}\label{le-w}
The increment $w_{q+1}=J_{q+1}-J_q$ satisfies the following estimates,
\begin{align}
\|w_{q+1}^p\|_{L^2}\lesssim&\ \delta_{q+1}^{\frac12}, \label{wp-l2}\\
\|w_{q+1}^c\|_{L^2}+\|w_{q+1}^t\|_{L^2}\lesssim&\ \ell^{-1}\mu^{-1}\delta_{q+1}^{\frac12}r^{\frac32}, \label{wc-l2}\\
\|w_{q+1}\|_{L^p}\lesssim&\ \delta_{q+1}^{\frac12}\ell^{-\frac12(1-\frac1p)}r^{\frac32-\frac3p}, \ \ p\geq 1, \label{w-p}\\
\|w_{q+1}^p\|_{W^{1,p}}+\|w_{q+1}^c\|_{W^{1,p}}+\|w_{q+1}^t\|_{W^{1,p}}\lesssim&\ \ell^{-2}\lambda_{q+1}r^{\frac32-\frac3p}, \ \ p\geq 1, \label{w-1p}\\
\|\partial_tw_{q+1}^p\|_{L^p}+\|\partial_tw_{q+1}^c\|_{L^p}\lesssim&\ \ell^{-2}\lambda_{q+1}\sigma\mu r^{\frac52-\frac3p}, \ \ p\geq 1, \label{wt-p}\\
\|w_{q+1}^p\|_{C^N_{x,t}}+\|w_{q+1}^c\|_{C^N_{x,t}}+\|w_{q+1}^t\|_{C^N_{x,t}}\lesssim&\ \lambda_{q+1}^{\frac{3+5N}{2}}, \label{w-cn}\\
\||\nabla|^Nw_{q+1}^p\|_{L^p}+\||\nabla|^Nw_{q+1}^c\|_{L^p}+& \||\nabla|^Nw_{q+1}^t\|_{L^p} \notag\\
\lesssim&\ \lambda_{q+1}^{N} r^{\frac32-\frac3p}, \ \ p\geq 1, \label{wn-p}\\
\||\nabla|^N J_{q+1}\|_{L^p}\lesssim&\ \lambda_{q+1}^{N+\frac32},  \ \ p\geq 1. \label{Jn-p}
\end{align}
\end{Lemma}
\pf
Recall that 
\begin{equation}\notag
\begin{split}
w_{q+1}^p=\nabla\times v_{q+1}^p
&\ = \sum_{\xi\in \Lambda}a_{\xi}\mathbb W_{\xi}+\lambda_{q+1}^{-1}\sum_{\xi\in \Lambda}\nabla(a_{\xi}\eta_{\xi})\times W_{\xi}\\
&\ =: \ \mathbb W^J+\mathbb W_{\epsilon,1}
\end{split}
\end{equation}
with 
\begin{equation}\label{eq-w-e1}
\mathbb W^J= \lambda_{q+1}v_{q+1}^p, \ \ \mathbb W_{\epsilon,1}=\lambda_{q+1}^{-1}\sum_{\xi\in \Lambda}\nabla(a_{\xi}\eta_{\xi})\times W_{\xi},
\end{equation}
and
\begin{equation}\notag
w_{q+1}^c=\nabla\times v_{q+1}^c
=\lambda_{q+1}^{-2}\nabla\times\left(\sum_{\xi\in \Lambda}\nabla(a_{\xi}\eta_{\xi})\times W_{\xi}\right)
= \lambda_{q+1}^{-1}\nabla\times\mathbb W_{\epsilon,1}.
\end{equation}
Note that $\mathbb W^J$ is the intermittent wave defined for the principle part of the velocity increment $u_{q+1}^p-u_q^p$ in \cite{BV}; while the temporal oscillation part $w_{q+1}^t$ is defined the same way as in \cite{BV}. Thus the estimates on $\mathbb W^J$ and $w_{q+1}^t$ can be adapted from \cite{BV}. Therefore, it is sufficient to estimate $\mathbb W_{\epsilon,1}$ and $\nabla\times\mathbb W_{\epsilon,1}$.

In addition, we notice that $\mathbb W_{\epsilon,1}=\lambda_{q+1} v_{q+1}^c$. It then follows from Lemma \ref{le-vq} that
\begin{equation}\label{est-We-sobolev}
\begin{split}
\|\mathbb W_{\epsilon,1}\|_{L^2}\leq&\ \lambda_{q+1}\|v_{q+1}^c\|_{L^2}\lesssim\delta_{q+1}^{\frac12},\\
\|\mathbb W_{\epsilon,1}\|_{L^p}\leq&\ \lambda_{q+1}\|v_{q+1}^c\|_{L^p}\lesssim \delta_{q+1}^{\frac12}\ell^{-\frac12(1-\frac1p)}\sigma r^{\frac52-\frac3p},\\
\|\mathbb W_{\epsilon,1}\|_{W^{1,p}}\leq&\ \lambda_{q+1}\| v_{q+1}^c\|_{W^{1,p}}\lesssim \ell^{-2}\lambda_{q+1}r^{\frac32-\frac3p},\\
\|\partial_t\mathbb W_{\epsilon,1}\|_{L^p}\leq&\ \lambda_{q+1}\|\partial_tv_{q+1}^c\|_{L^p}\lesssim \ell^{-2}\lambda_{q+1}\sigma\mu r^{\frac52-\frac3p},\\
\|\mathbb W_{\epsilon,1}\|_{C^N_{x,t}}\leq&\ \lambda_{q+1}\|v_{q+1}^c\|_{C^N_{x,t}}\lesssim \lambda_{q+1}^{\frac{3+5N}{2}}.\\
\end{split}
\end{equation}
The estimates of $\nabla\times\mathbb W_{\epsilon,1}$ and hence $w_{q+1}^c$ are carried out as follows. First, a direct computation leads to 
\begin{equation}\label{wc-re}
w_{q+1}^c=\lambda_{q+1}^{-2}\sum_{\xi\in\Lambda}\left(W_{\xi}\cdot \nabla\nabla(a_{\xi}\eta_{\xi})-W_{\xi}\cdot\Delta(a_{\xi}\eta_{\xi})-\nabla(a_{\xi}\eta_{\xi})\cdot\nabla W_{\xi}\right)
\end{equation}
where we used the fact that $\nabla\cdot W_{\xi}=0$. Thus, we have
\begin{equation}\notag
\begin{split}
\|w_{q+1}^c\|_{L^2}\lesssim &\ \lambda_{q+1}^{-2}\left(\|W_{\xi}\nabla\nabla(a_{\xi}\eta_{\xi})\|_{L^2}+\|\nabla(a_{\xi}\eta_{\xi})\nabla W_{\xi}\|_{L^2}\right)\\
\lesssim &\ \lambda_{q+1}^{-2}\left(\|\nabla\nabla(a_{\xi}\eta_{\xi})\|_{L^2}+\lambda_{q+1}\|\nabla(a_{\xi}\eta_{\xi})\|_{L^2}\right)\\
\lesssim &\  \lambda_{q+1}^{-2}\left(\|a_{\xi}\nabla^2\eta_{\xi}\|_{L^2}+\|\nabla a_{\xi}\nabla \eta_{\xi}\|_{L^2}+\|\nabla^2a_{\xi}\eta_{\xi}\|_{L^2}\right)\\
&+\lambda_{q+1}^{-1}\left(\|a_{\xi}\nabla \eta_{\xi}\|_{L^2}+\|\nabla a_{\xi}\eta_{\xi}\|_{L^2}\right).\\
\end{split}
\end{equation}
Following Lemma \ref{le-holder} for $L^2$ norm, we obtain
\begin{equation}\notag
\begin{split}
\|a_{\xi}\nabla^2\eta_{\xi}\|_{L^2}\lesssim&\ \|a_{\xi}\|_{L^2}\|\nabla^2\eta_{\xi}\|_{L^2}+\lambda_{q+1}^{-\frac12}\|a_{\xi}\|_{C^1_x}\|\nabla^2\eta_{\xi}\|_{L^2}\\
\lesssim&\ \delta_{q+1}^{\frac12}(\lambda_{q+1}\sigma r)^2+ \lambda_{q+1}^{-\frac12}\ell^{-1} (\lambda_{q+1}\sigma r)^2\\
\lesssim&\ \lambda_{q+1}^2\ell^{-1}\mu^{-1}\delta_{q+1}^{\frac12}r^{\frac32}
\end{split}
\end{equation}
due to (\ref{eta-lp}) and (\ref{a-l2}), and the choice of parameters (\ref{pp1}) and (\ref{pp2}). The other terms are treated in an analogous way,
\begin{equation}\notag
\begin{split}
\|\nabla a_{\xi}\nabla \eta_{\xi}\|_{L^2}\lesssim&\ \|\nabla a_{\xi}\|_{L^\infty}\|\nabla \eta_{\xi}\|_{L^2}\lesssim \ell^{-1}\lambda_{q+1}\sigma r\lesssim \lambda_{q+1}^2\ell^{-1}\mu^{-1}\delta_{q+1}^{\frac12}r^{\frac32};\\
\|\nabla^2 a_{\xi} \eta_{\xi}\|_{L^2}\lesssim&\ \|\nabla^2 a_{\xi}\|_{L^\infty}\| \eta_{\xi}\|_{L^2}\lesssim \ell^{-2}\lesssim \lambda_{q+1}^2\ell^{-1}\mu^{-1}\delta_{q+1}^{\frac12}r^{\frac32};\\
\|a_{\xi}\nabla\eta_{\xi}\|_{L^2}\lesssim&\ \|a_{\xi}\|_{L^2}\|\nabla\eta_{\xi}\|_{L^2}\lesssim \delta_{q+1}^{\frac12}\lambda_{q+1}\sigma r\lesssim \lambda_{q+1}\ell^{-1}\mu^{-1}\delta_{q+1}^{\frac12}r^{\frac32};\\
\|\nabla a_{\xi} \eta_{\xi}\|_{L^2}\lesssim&\ \|\nabla a_{\xi}\|_{L^\infty}\| \eta_{\xi}\|_{L^2}\lesssim \ell^{-1}\lesssim \lambda_{q+1}\ell^{-1}\mu^{-1}\delta_{q+1}^{\frac12}r^{\frac32}.
\end{split}
\end{equation}
Combining the estimates above yields 
\[\|w_{q+1}^c\|_{L^2}\lesssim \ell^{-1}\mu^{-1}\delta_{q+1}^{\frac12}r^{\frac32}\]
which concludes the proof of (\ref{wc-l2}).

Now we estimate the $L^p$ norm of $w_{q+1}^p$, $w_{q+1}^c$, and $w_{q+1}^t$. Again we recall that  $w_{q+1}^p=\lambda_{q+1}v_{q+1}^p+\lambda_{q+1}v_{q+1}^c$. The estimates (\ref{v-pp}) and (\ref{v-cp}) give immediately 
\begin{equation}\notag
\begin{split}
\|w_{q+1}^p\|_{L^p}\lesssim&\ \delta_{q+1}^{\frac12}\ell^{-\frac12(1-\frac1p)}r^{\frac32-\frac3p}+
\delta_{q+1}^{\frac12}\ell^{-\frac12(1-\frac1p)}\sigma r^{\frac52-\frac3p}\\
\lesssim&\ \delta_{q+1}^{\frac12}\ell^{-\frac12(1-\frac1p)}r^{\frac32-\frac3p}.
\end{split}
\end{equation}
In an analogous way of estimating $\|w_{q+1}^c\|_{L^2}$, we can obtain
\[\|w_{q+1}^c\|_{L^p}\lesssim \ell^{-1}\mu^{-1}\delta_{q+1}^{\frac12}r^{3-\frac3p}.\]
While we deal with $w_{q+1}^t$ as follows, by using (\ref{a-inf}) and (\ref{eta-lp})
\begin{equation}\notag
\begin{split}
\|w_{q+1}^t\|_{L^p}\lesssim &\ \mu^{-1}\sum_{\xi\in\Lambda_{\alpha}}\|a_{\xi}^2\eta_{\xi}^2\|_{L^p}\\
\lesssim &\ \mu^{-1}\sum_{\xi\in\Lambda_{\alpha}}\|a_{\xi}\|_{L^\infty}^2\|\eta_{\xi}\|_{L^{2p}}^2\\
\lesssim &\ \mu^{-1}\delta_{q+1}\ell^{-1}r^{3-\frac3p}.
\end{split}
\end{equation}
Combining the last three estimates yields
\[\|w_{q+1}\|_{L^p}\lesssim \delta_{q+1}^{\frac12}\ell^{-\frac12(1-\frac1p)}r^{\frac32-\frac3p}\]
which proves (\ref{w-p}).

Next we estimate $\|w_{q+1}^c\|_{W^{1,p}}$. It follows from (\ref{wc-re}), Lemma \ref{le-a}, Lemma \ref{le-w} and (\ref{pp1})-(\ref{pp2}) that
\begin{equation}\notag
\begin{split}
\|w_{q+1}^c\|_{W^{1,p}}\lesssim &\ \lambda_{q+1}^{-2}\left(\|W_{\xi}\nabla^2(a_{\xi}\eta_{\xi})\|_{W^{1,p}}+\|\nabla W_{\xi}\nabla (a_{\xi} \eta_{\xi})\|_{W^{1,p}}\right)\\
\lesssim &\ \lambda_{q+1}^{-2}\left(\|\nabla^3(a_{\xi}\eta_{\xi})\|_{L^p}+\lambda_{q+1}\|\nabla^2 (a_{\xi} \eta_{\xi})\|_{L^p}\right)\\
\lesssim &\ \lambda_{q+1}^{-2}\|a_{\xi}\|_{C^3}\left(\|\nabla^3\eta_{\xi}\|_{L^p}+\lambda_{q+1}\|\nabla^2  \eta_{\xi}\|_{L^p}\right)\\
\lesssim &\ \lambda_{q+1}^{-2} \ell^{-3}\left((\lambda_{q+1}\sigma r)^3+\lambda_{q+1}(\lambda_{q+1}\sigma r)^2\right)r^{\frac32-\frac3p}\\
\lesssim &\ \ell^{-2}\lambda_{q+1}r^{\frac32-\frac3p}.
\end{split}
\end{equation}
Thus, the proof of (\ref{w-1p}) is also complete.

To prove (\ref{wt-p}), we proceed to estimate $\|\partial_t w_{q+1}^c\|_{L^p}$,
\begin{equation}\notag
\begin{split}
\|\partial_t w_{q+1}^c\|_{L^p}\lesssim &\ \lambda_{q+1}^{-2}\sum_{\xi\in\Lambda}\left(\|\partial_t\nabla^2(a_{\xi}\eta_{\xi})\|_{L^p}+\lambda_{q+1}\|\partial_t\nabla(a_{\xi}\eta_{\xi})\|_{L^p}\right)\\
\lesssim &\ \lambda_{q+1}^{-2}\sum_{\xi\in\Lambda}\|a_{\xi}\|_{C^3}\left(\|\partial_t\nabla^2\eta_{\xi}\|_{L^p}+\|\nabla^2\eta_{\xi}\|_{L^p}\right)\\
 &+\lambda_{q+1}^{-1}\sum_{\xi\in\Lambda}\|a_{\xi}\|_{C^2}\left(\|\partial_t\nabla\eta_{\xi}\|_{L^p}+\|\nabla\eta_{\xi}\|_{L^p}\right)\\
\lesssim &\ \lambda_{q+1}^{-2} \ell^{-3}\left(\lambda_{q+1}\sigma r\mu+1\right)(\lambda_{q+1}\sigma r)^2r^{\frac32-\frac3p}\\
&+\lambda_{q+1}^{-1} \ell^{-2}\left(\lambda_{q+1}\sigma r\mu+1\right)(\lambda_{q+1}\sigma r)r^{\frac32-\frac3p}\\
\lesssim &\ \ell^{-2}\lambda_{q+1}\sigma\mu r^{\frac52-\frac3p}.
\end{split}
\end{equation}
In the end, we estimate $\|w_{q+1}^c\|_{C^N_{x,t}}$,
\begin{equation}\notag
\begin{split}
\|w_{q+1}^c\|_{C^N_{x,t}}\lesssim &\ \lambda_{q+1}^{-2}\sum_{\xi\in\Lambda}\left(\|\nabla^2(a_{\xi}\eta_{\xi})\|_{C^N_{x,t}}+\lambda_{q+1}\|\nabla(a_{\xi}\eta_{\xi})\|_{C^N_{x,t}}\right)\\
\lesssim &\ \lambda_{q+1}^{-2}\sum_{\xi\in\Lambda}\left(\|a_{\xi}\|_{C^{N}_{x,t}}\|\nabla^2\eta_{\xi}\|_{C^{N}_{x,t}}+\|\nabla a_{\xi}\|_{C^{N}_{x,t}}\|\nabla\eta_{\xi}\|_{C^{N}_{x,t}}+\|\nabla^2a_{\xi}\|_{C^{N}_{x,t}}\|\eta_{\xi}\|_{C^{N}_{x,t}}\right)\\
&+\lambda_{q+1}^{-1}\sum_{\xi\in\Lambda}\left(\|a_{\xi}\|_{C^{N}_{x,t}}\|\nabla\eta_{\xi}\|_{C^{N}_{x,t}}+\|\nabla a_{\xi}\|_{C^{N}_{x,t}}\|\eta_{\xi}\|_{C^{N}_{x,t}}\right)\\
\lesssim &\ \lambda_{q+1}^{-2}\left(\ell^{-N}(\lambda_{q+1}\sigma r)^2 +\ell^{-N-1}\lambda_{q+1}\sigma r+\ell^{-N-2}\right)(\lambda_{q+1}\sigma r\mu)^Nr^{\frac32}\\
&+\lambda_{q+1}^{-1}\left(\ell^{-N}\lambda_{q+1}\sigma r+\ell^{-N-1}\right)(\lambda_{q+1}\sigma r\mu)^Nr^{\frac32}\\
\lesssim &\ \lambda_{q+1}^{\frac32+\frac{5N}2}.
\end{split}
\end{equation}
It completes the proof of the inequality (\ref{w-cn}).

Inequality (\ref{wn-p}) can be obtained analogously as (\ref{w-1p}); while (\ref{Jn-p}) is implied by (\ref{wn-p}) and (\ref{pp2}).

\cbdu

\medskip

\subsection{Stress tensor $R_{q+1}$ and its estimate}
\label{sec-emhd-rq1}

Based on the construction of the increments $v_{q+1}=B_{q+1}-B_q$ and $w_{q+1}=J_{q+1}-J_q$, and hence $B_{q+1}=B_{q}+v_{q+1}$ and $J_{q+1}=J_{q}+w_{q+1}$, we will derive the new stress tensor $R_{q+1}$ such that $(B_{q+1}, J_{q+1}, R_{q+1})$ satisfies
\begin{equation}\label{eq-emhd-R-q1}
\partial_tJ_{q+1}+\nabla\times\nabla\times (J_{q+1}\times B_{q+1})=\Delta J_{q+1}+\nabla\cdot  R_{q+1},
\end{equation}
i.e. equation (\ref{eq-div-q-alt}) (and (\ref{eq-div-q}) as well) at the level $q+1$. 

Substituting $B_{q+1}=B_{q}+v_{q+1}$ and $J_{q+1}=J_{q}+w_{q+1}$ in the left hand side of (\ref{eq-emhd-R-q1}), we have
\begin{equation}\label{eq-emhd-R-q2}
\begin{split}
&\partial_tJ_{q+1}+\nabla\times\nabla\times (J_{q+1}\times B_{q+1})\\
=& \partial_tJ_{q}+\partial_t w_{q+1}+\nabla\times\nabla\times ((J_q+w_{q+1})\times (B_q+v_{q+1}))\\
=&  \partial_tJ_{q}+\partial_t w_{q+1}+\nabla\times\nabla\times (J_q\times B_q)
+\nabla\times\nabla\times (J_q\times v_{q+1})\\
&+\nabla\times\nabla\times (w_{q+1}\times B_q)+\nabla\times\nabla\times (w_{q+1}\times v_{q+1}).
\end{split}
\end{equation}
Recall that $(B_{q}, J_{q}, R_{q})$ satisfies equation (\ref{eq-div-q-alt}), i.e. 
\begin{equation}\notag
\partial_tJ_{q}+\nabla\times\nabla\times (J_{q}\times B_{q})=\Delta J_{q}+\nabla\cdot  R_{q}.
\end{equation}
Since $\Delta J_q=\Delta J_{q+1}-\Delta w_{q+1}$, we have
\begin{equation}\notag
\partial_tJ_{q}+\nabla\times\nabla\times (J_{q}\times B_{q})=\Delta J_{q+1}-\Delta w_{q+1}+\nabla\cdot  R_{q}.
\end{equation}
Therefore, substituting $\partial_tJ_{q}+\nabla\times\nabla\times (J_{q}\times B_{q})$ by $\Delta J_{q+1}-\Delta w_{q+1}+\nabla\cdot  R_{q}$ in the right hand side of (\ref{eq-emhd-R-q2}) leads to

\begin{equation}\label{eq-emhd-R-q3}
\begin{split}
&\partial_tJ_{q+1}+\nabla\times\nabla\times (J_{q+1}\times B_{q+1})\\
=&\Delta J_{q+1} +\nabla\cdot  R_{q}+ \partial_t w_{q+1}-\Delta w_{q+1}
+\nabla\times\nabla\times (J_q\times v_{q+1})\\
&+\nabla\times\nabla\times (w_{q+1}\times B_q)+\nabla\times\nabla\times (w_{q+1}\times v_{q+1}).
\end{split}
\end{equation}
We now observe from (\ref{eq-emhd-R-q3}) that, if we choose $R_{q+1}$ such that 
\begin{equation}\label{eq-emhd-R-q4}
\begin{split}
\nabla\cdot  R_{q+1}=&\nabla\cdot  R_{q}+ \partial_t w_{q+1}-\Delta w_{q+1}
+\nabla\times\nabla\times (J_q\times v_{q+1})\\
&+\nabla\times\nabla\times (w_{q+1}\times B_q)+\nabla\times\nabla\times (w_{q+1}\times v_{q+1}),
\end{split}
\end{equation}
then $(B_{q+1}, J_{q+1}, R_{q+1})$ satisfies (\ref{eq-emhd-R-q1}). 

Next, we rewrite the double curl terms with $\nabla\times \nabla\times $ in (\ref{eq-emhd-R-q4}) into divergence form of $\nabla\cdot$. 
Applying a few vector identities in Section \ref{sec-vec}, we have
\begin{equation}\notag
\begin{split}
&\nabla\times \nabla\times (J_q\times v_{q+1})\\
=&\ \nabla\cdot[v_{q+1}\otimes (\nabla\times J_q)-(\nabla\times J_q)\otimes v_{q+1}]\\
& +\nabla\cdot[(\nabla\times v_{q+1})\otimes  J_q-J_q\otimes (\nabla\times v_{q+1})]\\
& - 2\nabla\times (J_q\nabla v_{q+1}),
\end{split}
\end{equation}
\begin{equation}\notag
\begin{split}
&\nabla\times \nabla\times (w_{q+1}\times B_{q})\\
=&\ \nabla\cdot[B_{q}\otimes (\nabla\times w_{q+1})-(\nabla\times w_{q+1})\otimes B_{q}]\\
& +\nabla\cdot[(\nabla\times B_{q})\otimes  w_{q+1}-w_{q+1}\otimes (\nabla\times B_{q})]\\
& - 2\nabla\times (w_{q+1}\nabla B_{q}),
\end{split}
\end{equation}
\begin{equation}\notag
\begin{split}
&\nabla\times \nabla\times (w_{q+1}\times v_{q+1})\\
=&\ \nabla\cdot[v_{q+1}\otimes (\nabla\times w_{q+1})-(\nabla\times w_{q+1})\otimes v_{q+1}]\\
& +\nabla\cdot[(\nabla\times v_{q+1})\otimes  w_{q+1}-w_{q+1}\otimes (\nabla\times v_{q+1})]\\
& - 2\nabla\times (w_{q+1}\nabla v_{q+1}).
\end{split}
\end{equation}
Combining the last three equations with (\ref{eq-emhd-R-q4}), 
we obtain
\begin{equation}\label{eq-emhd-R-q5}
\begin{split}
\nabla\cdot  R_{q+1}=&\ \partial_t w_{q+1}-\Delta w_{q+1}\\
&+\nabla\cdot[v_{q+1}\otimes (\nabla\times J_q)-(\nabla\times J_q)\otimes v_{q+1}]\\
& +\nabla\cdot[(\nabla\times v_{q+1})\otimes  J_q-J_q\otimes (\nabla\times v_{q+1})]\\
&+ \nabla\cdot[B_{q}\otimes (\nabla\times w_{q+1})-(\nabla\times w_{q+1})\otimes B_{q}]\\
& +\nabla\cdot[(\nabla\times B_{q})\otimes  w_{q+1}-w_{q+1}\otimes (\nabla\times B_{q})]\\
&+ \nabla\cdot[v_{q+1}\otimes (\nabla\times w_{q+1})-(\nabla\times w_{q+1})\otimes v_{q+1}]\\
& +\nabla\cdot[(\nabla\times v_{q+1})\otimes  w_{q+1}-w_{q+1}\otimes (\nabla\times v_{q+1})]\\
&-2\nabla\times (J_q\nabla v_{q+1})- 2\nabla\times (w_{q+1}\nabla B_{q})- 2\nabla\times (w_{q+1}\nabla v_{q+1}) \\
&+\nabla\cdot  R_{q}.
\end{split}
\end{equation}
We further reordering the terms and classify them on the right hand side of (\ref{eq-emhd-R-q5}) into linear, correction and oscillation terms:
\begin{equation}\label{eq-emhd-R-q6}
\begin{split}
&\nabla\cdot R_{q+1}\\
=&\ \{\nabla\cdot[\mathcal R(\partial _t w_{q+1}^p+\partial _t w_{q+1}^c-\Delta w_{q+1})]\\
&+\nabla\cdot[v_{q+1}\otimes(\nabla\times J_q)-(\nabla\times J_q)\otimes v_{q+1}
+(\nabla\times v_{q+1})\otimes J_q-J_q\otimes (\nabla\times v_{q+1})]\\
&+\nabla\cdot[B_q\otimes (\nabla\times w_{q+1})-(\nabla\times w_{q+1})\otimes B_q
+(\nabla\times B_q)\otimes w_{q+1}-w_{q+1}\otimes (\nabla\times B_q)]\\
&-\nabla\cdot[\mathcal R(2\nabla\times (J_q\nabla v_{q+1})+ 2\nabla\times (w_{q+1}\nabla B_{q}))]\}\\
&+\{\nabla\cdot[(v_{q+1}^c+v_{q+1}^t)\otimes (\nabla\times w_{q+1})+ v_{q+1}^p\otimes(\nabla\times(w_{q+1}^c+w_{q+1}^t))]\\
&-\nabla\cdot[(\nabla\times(w_{q+1}^c+w_{q+1}^t))\otimes v_{q+1}+(\nabla\times w_{q+1}^p)\otimes(v_{q+1}^c+v_{q+1}^t)]\\
&+\nabla\cdot[ (\nabla\times(v_{q+1}^c+v_{q+1}^t))\otimes w_{q+1}+(\nabla\times v_{q+1}^p)\otimes(w_{q+1}^c+w_{q+1}^t)]\\
&-\nabla\cdot[ (w_{q+1}^c+w_{q+1}^t)\otimes (\nabla\times v_{q+1})-w_{q+1}^p\otimes(\nabla\times(v_{q+1}^c+v_{q+1}^t))]\\
&-\nabla\cdot[\mathcal R\nabla\times (2w_{q+1}\nabla v_{q+1}^c)]-2\nabla\times (w_{q+1}\nabla v_{q+1}^t)\\
&-\nabla\cdot[ \mathcal R\nabla\times(2(w_{q+1}^c+w_{q+1}^t)\nabla v_{q+1}^p)]\}\\
&+\{\nabla\cdot[v_{q+1}^p\otimes (\nabla\times w_{q+1}^p)-(\nabla\times w_{q+1}^p)\otimes v_{q+1}^p]\\
&+\nabla\cdot[(\nabla\times v_{q+1}^p)\otimes  w_{q+1}^p-w_{q+1}^p\otimes (\nabla\times v_{q+1}^p)]\\
&-2\nabla\times(w_{q+1}^p\nabla v_{q+1}^p)
+(\nabla\cdot R_q+\partial_tw_{q+1}^t)\} \\
=:&\  \nabla\cdot R_{\mathrm{linear}}+\nabla\cdot R_{\mathrm{corrector}}+\nabla\cdot R_{\mathrm{oscillation}}.
\end{split}
\end{equation}
On the right hand side of the equation (\ref{eq-emhd-R-q6}), the first four lines correspond to linear terms, the middle six lines correspond to correction terms, and the last three lines correspond to oscillation terms. 

Before diving into the estimates of the new stress tensor, we will first rearrange $\nabla\cdot R_{\mathrm{corrector}}$ and $\nabla\cdot R_{\mathrm{oscillation}}$. 

\begin{Lemma}\label{le-rewrite-r-corr}
The term $\nabla\cdot R_{\mathrm{corrector}}$ can be reformulated as
\begin{equation}\label{eq-r-corrector2}
\begin{split}
\nabla\cdot R_{\mathrm{corrector}}
=&  \nabla\cdot (R_{\{\mathrm{cor},2\}}+R_{\{\mathrm{cor},3\}})+\nabla\times\nabla\times \tilde M_{q+1}+\nabla\tilde p_{q+1,1}\\
\end{split}
\end{equation}
with $R_{\{\mathrm{cor},2\}}$, $R_{\{\mathrm{cor},3\}}$ and $\tilde M_{q+1}$ defined as follows
\begin{equation}\label{eq-r-corr-2}
\begin{split}
R_{\{\mathrm{cor},2\}}=&\ \nabla\cdot[v_{q+1}^c\otimes (\nabla\times w_{q+1})-(\nabla\times w_{q+1}^p)\otimes v_{q+1}^c]\\
&+\nabla\cdot[(\nabla\times v_{q+1}^c)\otimes w_{q+1}-w_{q+1}\otimes (\nabla\times v_{q+1}^c) ]\\
&-\nabla\cdot[\mathcal R\nabla\times (2w_{q+1}\nabla v_{q+1}^c)],
\end{split}
\end{equation} 
\begin{equation}\label{eq-r-corr-3}
\begin{split}
R_{\{\mathrm{cor},3\}}=&\ \nabla\cdot[(\nabla\times v_{q+1}^p)\otimes (w_{q+1}^c+w_{q+1}^t)]\\
&-\nabla\cdot[(w_{q+1}^c+w_{q+1}^t)\otimes (\nabla\times v_{q+1}^p)]\\
&-\nabla\cdot[ \mathcal R\nabla\times(2(w_{q+1}^c+w_{q+1}^t)\nabla v_{q+1}^p)],
\end{split}
\end{equation} 
\begin{equation}\label{eq-m-q11}
\tilde M_{q+1}=w_{q+1}\times v_{q+1}^t,
\end{equation}
while $\tilde p_{q+1,1}$ is a dummy pressure term to make $R_{\{\mathrm{cor},2\}}$ and $R_{\{\mathrm{cor},3\}}$ traceless.
\end{Lemma}
\pf
Upon the choice of parameters (\ref{pp1})-(\ref{pp2}), the upper bound of $\|v_{q+1}^t\|_{L^p}$ is larger than that of $\|v_{q+1}^c\|_{L^p}$, as in Lemma \ref{le-vq}. Thus we have to handle the terms involving $\|v_{q+1}^t\|_{L^p}$ in $R_{\mathrm{corrector}}$ in a more delicate way. In order to do so, we rearrange $\nabla\cdot R_{\mathrm{corrector}}$ as,
\begin{equation}\notag
\begin{split}
&\nabla\cdot R_{\mathrm{corrector}}\\
=&\ \nabla\cdot[(v_{q+1}^c+v_{q+1}^t)\otimes (\nabla\times w_{q+1})+ v_{q+1}^p\otimes(\nabla\times(w_{q+1}^c+w_{q+1}^t))]\\
&-\nabla\cdot[(\nabla\times(w_{q+1}^c+w_{q+1}^t))\otimes v_{q+1}+(\nabla\times w_{q+1}^p)\otimes(v_{q+1}^c+v_{q+1}^t)]\\
&+\nabla\cdot[ (\nabla\times(v_{q+1}^c+v_{q+1}^t))\otimes w_{q+1}+(\nabla\times v_{q+1}^p)\otimes(w_{q+1}^c+w_{q+1}^t)]\\
&-\nabla\cdot[ (w_{q+1}^c+w_{q+1}^t)\otimes (\nabla\times v_{q+1})-w_{q+1}^p\otimes(\nabla\times(v_{q+1}^c+v_{q+1}^t))]\\
&-\nabla\cdot[\mathcal R\nabla\times (2w_{q+1}\nabla v_{q+1}^c)]-2\nabla\times (w_{q+1}\nabla v_{q+1}^t)\\
&-\nabla\cdot[ \mathcal R\nabla\times(2(w_{q+1}^c+w_{q+1}^t)\nabla v_{q+1}^p)]\}\\
=&\ \left\{\nabla\cdot[v_{q+1}^t\otimes \nabla\times w_{q+1}-(\nabla\times w_{q+1}^p)\otimes v_{q+1}^t-(\nabla\times (w_{q+1}^c+w_{q+1}^t))\otimes v_{q+1}^t]\right.\\
&+\nabla\cdot[(\nabla\times v_{q+1}^t)\otimes w_{q+1}-w_{q+1}^p\otimes (\nabla\times v_{q+1}^t)-(w_{q+1}^c+w_{q+1}^t)\otimes(\nabla\times v_{q+1}^t)]\\
&\left.-2\nabla\times (w_{q+1}\nabla v_{q+1}^t)\right\}\\
&+\left\{\nabla\cdot[v_{q+1}^c\otimes (\nabla\times w_{q+1})-(\nabla\times w_{q+1}^p)\otimes v_{q+1}^c]\right.\\
&+\left.\nabla\cdot[(\nabla\times v_{q+1}^c)\otimes w_{q+1}-w_{q+1}\otimes (\nabla\times v_{q+1}^c) ]-\nabla\cdot[\mathcal R\nabla\times (2w_{q+1}\nabla v_{q+1}^c)] \right\}\\
&+\left\{\nabla\cdot[(\nabla\times v_{q+1}^p)\otimes (w_{q+1}^c+w_{q+1}^t)-(w_{q+1}^c+w_{q+1}^t)\otimes (\nabla\times v_{q+1}^p)]\right.\\
&\left.-\nabla\cdot[ \mathcal R\nabla\times(2(w_{q+1}^c+w_{q+1}^t)\nabla v_{q+1}^p)]\right\}\\
=&: \nabla\cdot R_{\{\mathrm{cor},1\}}+\nabla\cdot R_{\{\mathrm{cor},2\}}+\nabla\cdot R_{\{\mathrm{cor},3\}}.
\end{split}
\end{equation}
We notice that only $R_{\{\mathrm{cor},1\}}$ involves with $v_{q+1}^t$. We can further rewrite $\nabla\cdot R_{\{\mathrm{cor},1\}}$ into
\begin{equation}\notag
\begin{split}
 &\nabla\cdot R_{\{\mathrm{cor},1\}}\\
 =&\ \left\{\nabla\cdot[v_{q+1}^t\otimes \nabla\times w_{q+1}-(\nabla\times w_{q+1})\otimes v_{q+1}^t]\right.\\
&+\nabla\cdot[(\nabla\times v_{q+1}^t)\otimes w_{q+1}-w_{q+1}\otimes (\nabla\times v_{q+1}^t)]\\
&\left.-2\nabla\times (w_{q+1}\nabla v_{q+1}^t)\right\}\\
=&\nabla\times\nabla\times (w_{q+1}\times v_{q+1}^t).\\
\end{split}
\end{equation}
Denote 
\begin{equation}\notag
\tilde M_{q+1}=w_{q+1}\times v_{q+1}^t.
\end{equation}
It follows that
\begin{equation}\notag
 \nabla\cdot R_{\{\mathrm{cor},1\}}=\nabla\times\nabla\times \tilde M_{q+1}.
\end{equation}
It completes the proof of the lemma.
\cbdu

\begin{Lemma}\label{le-rewrite-r-osci}
The oscillation part $\nabla\cdot R_{\mathrm{oscillation}}$ of the stress tensor can be written as 
\begin{equation}\label{eq-r-osc-rewrite}
\begin{split}
\nabla\cdot R_{\mathrm{oscillation}}
=&  \nabla\cdot (R_{\{\mathrm{osc},1\}}+R_{\{\mathrm{osc},2\}})+\nabla \tilde p_{q+1,2}
\end{split}
\end{equation}
with $R_{\{\mathrm{osc},1\}}$ and $R_{\{\mathrm{osc},2\}}$ defined as follows
\begin{equation}\label{eq-r-osc-1}
R_{\{\mathrm{osc},1\}}=\lambda_{q+1}v_{q+1}^p\otimes \lambda_{q+1}v_{q+1}^p-w_{q+1}^p\otimes w_{q+1}^p+ R_q+\mathcal R\partial_tw_{q+1}^t,
\end{equation} 
\begin{equation}\label{eq-r-osc-2}
\begin{split}
R_{\{\mathrm{osc},2\}}=&\ 2\nabla\cdot\left(w_{q+1}^p\otimes \mathbb W_{\epsilon,1}\right)
+\nabla\cdot\left(w_{q+1}^c\otimes v_{q+1}^p-v_{q+1}^p\otimes w_{q+1}^c\right)\\
&-\nabla\cdot[\mathcal R\nabla\times(2\mathbb W_{\epsilon,1}\nabla v_{q+1}^p]
-\nabla\cdot\left(\mathbb W_{\epsilon,1}\otimes \mathbb W_{\epsilon,1}\right),
\end{split}
\end{equation} 
while $\mathbb W_{\epsilon,1}$ is defined in (\ref{eq-w-e1}) and $\tilde p_{q+1,2}$ is a pressure term to make $R_{\{\mathrm{osc},1\}}$ and $R_{\{\mathrm{osc},2\}}$ traceless.
\end{Lemma}

\pf
In fact the oscillation terms, see (\ref{eq-emhd-R-q6}), can be written as 
\begin{equation}\label{R-os}
\begin{split}
&\nabla\cdot R_{\mathrm{oscillation}}\\
=&\nabla\times\nabla\times (w_{q+1}^p\times v_{q+1}^p)+\nabla\cdot R_q+\partial_tw_{q+1}^t\\
=&\nabla\times\nabla\times ((\nabla\times v_{q+1}^p)\times v_{q+1}^p)+\nabla\cdot R_q+\partial_tw_{q+1}^t\\
=&\nabla\cdot(\lambda_{q+1}v_{q+1}^p\otimes \lambda_{q+1}v_{q+1}^p)-\nabla\cdot (w_{q+1}^p\otimes w_{q+1}^p)+\nabla\cdot R_q+\partial_tw_{q+1}^t\\ 
&+\nabla\times\nabla\times ((\nabla\times v_{q+1}^p)\times v_{q+1}^p)-\nabla\cdot(\lambda_{q+1}v_{q+1}^p\otimes \lambda_{q+1}v_{q+1}^p)\\
&+\nabla\cdot (w_{q+1}^p\otimes w_{q+1}^p).
\end{split}
\end{equation}
We denote the first three terms in (\ref{R-os}) by $\nabla\cdot R_{\{\mathrm{osc},1\}}$, i.e.
\[R_{\{\mathrm{osc},1\}}=\lambda_{q+1}v_{q+1}^p\otimes \lambda_{q+1}v_{q+1}^p- w_{q+1}^p\otimes w_{q+1}^p+ R_q+\mathcal R\partial_tw_{q+1}^t.\]
Regarding the last three terms of $\nabla\cdot R_{\mathrm{oscillation}}$ in (\ref{R-os}), we need more effort to deal with them.
We first recall that
\[w_{q+1}^p=\nabla\times v_{q+1}^p=\lambda_{q+1}v_{q+1}^p+\lambda_{q+1}^{-1}\sum_{\xi\in\Lambda}\nabla(a_{\xi}\eta_{\xi})\times W_{\xi}=\lambda_{q+1}v_{q+1}^p+\mathbb W_{\epsilon,1}.\]
Thus, we have 
\[\lambda_{q+1}v_{q+1}^p=w_{q+1}^p-\mathbb W_{\epsilon,1}.\]
On the other hand, we notice that
\begin{equation}\notag
\begin{split}
\nabla\times w_{q+1}^p=&\lambda_{q+1}\sum_{\xi\in \Lambda}a_{\xi}\mathbb W_{\xi}+\sum_{\xi\in\Lambda}\nabla(a_{\xi}\eta_{\xi})\times W_{\xi}+\nabla\times \mathbb W_{\epsilon,1}\\
=&\lambda_{q+1}\left(w_{q+1}^p-\mathbb W_{\epsilon,1}\right)+\lambda_{q+1}\mathbb W_{\epsilon,1}+\nabla\times \mathbb W_{\epsilon,1}\\
=&\lambda_{q+1}w_{q+1}^p+\nabla\times \mathbb W_{\epsilon,1}.\\
\end{split}
\end{equation}
Therefore, a straightforward computation leads to 
\begin{equation}\notag
\begin{split}
&\nabla\times\nabla\times ((\nabla\times v_{q+1}^p)\times v_{q+1}^p)-\nabla\cdot(\lambda_{q+1}v_{q+1}^p\otimes \lambda_{q+1}v_{q+1}^p)+\nabla\cdot (w_{q+1}^p\otimes w_{q+1}^p)\\
=&\ \nabla\times\nabla\times (w_{q+1}^p\times v_{q+1}^p)-\nabla\cdot\left((w_{q+1}^p-\mathbb W_{\epsilon,1})\otimes (w_{q+1}^p-\mathbb W_{\epsilon,1})\right)\\
&+\nabla\cdot (w_{q+1}^p\otimes w_{q+1}^p)\\
=&\ \nabla\cdot\left(v_{q+1}^p\otimes(\nabla\times w_{q+1}^p)-(\nabla\times w_{q+1}^p)\otimes v_{q+1}^p\right)-2\nabla\times (w_{q+1}^p\nabla v_{q+1}^p)\\
&+\nabla\cdot\left(\mathbb W_{\epsilon,1}\otimes w_{q+1}^p\right)
+\nabla\cdot\left(w_{q+1}^p\otimes \mathbb W_{\epsilon,1}\right)-\nabla\cdot\left(\mathbb W_{\epsilon,1}\otimes \mathbb W_{\epsilon,1}\right)\\
=&\ \nabla\cdot\left(v_{q+1}^p\otimes (\lambda_{q+1}w_{q+1}^p-\nabla\times \mathbb W_{\epsilon,1})-(\lambda_{q+1}w_{q+1}^p-\nabla\times \mathbb W_{\epsilon,1})\otimes v_{q+1}^p\right)\\
&-2\nabla\times\left(\mathbb W_{\epsilon,1}\nabla v_{q+1}^p\right)
+\nabla\cdot\left(\mathbb W_{\epsilon,1}\otimes w_{q+1}^p\right)\\
&+\nabla\cdot\left(w_{q+1}^p\otimes \mathbb W_{\epsilon,1}\right)-\nabla\cdot\left(\mathbb W_{\epsilon,1}\otimes \mathbb W_{\epsilon,1}\right)\\
=&\ \nabla\cdot\left(\lambda_{q+1} v_{q+1}^p\otimes w_{q+1}^p-w_{q+1}^p\otimes \lambda_{q+1} v_{q+1}^p\right)\\
&-\nabla\cdot\left(v_{q+1}^p\otimes (\nabla\times \mathbb W_{\epsilon,1})-(\nabla\times \mathbb W_{\epsilon,1})\otimes v_{q+1}^p\right)\\
&-2\nabla\times\left(\mathbb W_{\epsilon,1}\nabla v_{q+1}^p\right)
+\nabla\cdot\left(\mathbb W_{\epsilon,1}\otimes w_{q+1}^p\right)\\
&+\nabla\cdot\left(w_{q+1}^p\otimes \mathbb W_{\epsilon,1}\right)-\nabla\cdot\left(\mathbb W_{\epsilon,1}\otimes \mathbb W_{\epsilon,1}\right)\\
=&\ \nabla\cdot\left((w_{q+1}^p-\mathbb W_{\epsilon,1})\otimes w_{q+1}^p-w_{q+1}^p\otimes (w_{q+1}^p-\mathbb W_{\epsilon,1})\right)\\
&-\nabla\cdot\left(v_{q+1}^p\otimes (\nabla\times \mathbb W_{\epsilon,1})-(\nabla\times \mathbb W_{\epsilon,1})\otimes v_{q+1}^p\right)\\
&-2\nabla\times\left(\mathbb W_{\epsilon,1}\nabla v_{q+1}^p\right)
+\nabla\cdot\left(\mathbb W_{\epsilon,1}\otimes w_{q+1}^p\right)\\
&+\nabla\cdot\left(w_{q+1}^p\otimes \mathbb W_{\epsilon,1}\right)-\nabla\cdot\left(\mathbb W_{\epsilon,1}\otimes \mathbb W_{\epsilon,1}\right).
\end{split}
\end{equation}
In view of the facts that $\nabla\cdot w_{q+1}^p=0$ and $\nabla\times \mathbb W_{\epsilon,1}=w_{q+1}^c$, we obtain that by continuing with the last equation
\begin{equation}\label{eq-d-q1-2}
\begin{split}
&\nabla\times\nabla\times ((\nabla\times v_{q+1}^p)\times v_{q+1}^p)-\nabla\cdot(\lambda_{q+1}v_{q+1}^p\otimes \lambda_{q+1}v_{q+1}^p)+\nabla\cdot (w_{q+1}^p\otimes w_{q+1}^p)\\
=&\ 2\nabla\cdot\left(w_{q+1}^p\otimes \mathbb W_{\epsilon,1}\right)
+\nabla\cdot\left(w_{q+1}^c\otimes v_{q+1}^p-v_{q+1}^p\otimes w_{q+1}^c\right)\\
&-\nabla\cdot[\mathcal R\nabla\times(2\mathbb W_{\epsilon,1}\nabla v_{q+1}^p]
-\nabla\cdot\left(\mathbb W_{\epsilon,1}\otimes \mathbb W_{\epsilon,1}\right).
\end{split}
\end{equation}
Thus, the conclusion of the lemma follows from (\ref{R-os}) and (\ref{eq-d-q1-2}).

\cbdu

We are now ready to accomplish the estimate of the new stress tensor $R_{q+1}$.

\begin{Lemma}\label{le-R}
Consider equation (\ref{eq-emhd-R-q1}) with $R_{q+1}$ defined by (\ref{eq-emhd-R-q5}) (and equivalently (\ref{eq-emhd-R-q6})).
There exists another
 traceless symmetric tensor $\tilde R_{q+1}$  
 and a scalar pressure function $\tilde p_{q+1}$ such that $\nabla\cdot R_{q+1}$ can be written as
\begin{equation}\label{eq-r-q1-final1}
\nabla\cdot R_{q+1}=\nabla\cdot \tilde R_{q+1}+\nabla\times\nabla\times \tilde M_{q+1}+\nabla\tilde p_{q+1}.
\end{equation}
In addition, there exists $p>1$ sufficiently close to 1, and a sufficiently small $\varepsilon_R>0$ independent of $q$ such that 
\begin{equation}\label{est-RM}
\|\tilde R_{q+1}\|_{L^p}+\|\tilde M_{q+1}\|_{L^p}\lesssim \lambda_{q+1}^{-2\varepsilon_R}\delta_{q+2}
\end{equation}
holds for some implicit constant which depends on $p$ and $\varepsilon_R$.
\end{Lemma}

\pf
Recall from (\ref{eq-emhd-R-q6}) that
\[ \nabla\cdot R_{q+1}=\nabla\cdot R_{\mathrm{linear}}+\nabla\cdot R_{\mathrm{corrector}}+\nabla\cdot R_{\mathrm{oscillation}}.\]
Denote 
\begin{equation}\label{eq-reformed-r-q1}
\begin{split}
\tilde R_{q+1}= &\ R_{\mathrm{linear}}+R_{\{\mathrm{cor},2\}}+R_{\{\mathrm{cor},3\}}+R_{\{\mathrm{osc},1\}}+R_{\{\mathrm{osc},2\}},\\
\tilde p_{q+1}=&\ \tilde p_{q+1,1}+\tilde p_{q+1,2}.
\end{split}
\end{equation}
Therefore, (\ref{eq-r-q1-final1}) follows from Lemma \ref{le-rewrite-r-corr}, Lemma \ref{le-rewrite-r-osci}, and (\ref{eq-reformed-r-q1}).

The estimate of $\tilde R_{q+1}$ will be established in Lemma \ref{le-linear}, Lemma \ref{le-corrector}, and Lemma \ref{le-oscill} below, and the estimate of $\tilde M_{q+1}$ will be obtained in Lemma \ref{le-corrector}. 
 Thus, estimate (\ref{est-RM}) follows from these lemmas.

\cbdu

\medskip

\subsubsection{Linear terms.}
The estimates of the linear terms are relatively easy.
\begin{Lemma}\label{le-linear}
For $p>1$ sufficiently close to $1$, $R_{\mathrm{linear}}$ satisfies
\begin{equation}\notag
\|R_{\mathrm{linear}}\|_{L^p}\lesssim \lambda_{q+1}^{-2\varepsilon_R}\delta_{q+2}. 
\end{equation}
\end{Lemma}
\pf
It follows from Lemma \ref{le-anti} and (\ref{w-1p}) that,
\begin{equation}\notag
\|\mathcal R\Delta w_{q+1}\|_{L^p}\lesssim \|w_{q+1}\|_{W^{1,p}}\lesssim \ell^{-2}\lambda_{q+1}r^{\frac32-\frac3p};
\end{equation}
while Lemma \ref{le-anti} and (\ref{wt-p}) together give
\begin{equation}\notag
\begin{split}
\|\mathcal R(\partial_t (w_{q+1}^p+w_{q+1}^c))\|_{L^p}=&\ \|\mathcal R(\partial_t \nabla\times (v_{q+1}^p+v_{q+1}^c))\|_{L^p}\\
=&\ \lambda_{q+1}^{-1} \|\mathcal R\partial_t \nabla\times\nabla\times v_{q+1}^p\|_{L^p}\\
=&\ \lambda_{q+1}^{-1} \|\mathcal R\partial_t \nabla\times w_{q+1}^p\|_{L^p}\\
\lesssim &\ \lambda_{q+1}^{-1} \|\partial_t  w_{q+1}^p\|_{L^p}\\
\lesssim &\ \ell^{-2}\sigma\mu r^{\frac52-\frac3p}.
\end{split}
\end{equation}
We have, by (\ref{Jn-p}) and (\ref{v-pp})-(\ref{v-tp}),
\begin{equation}\notag
\begin{split}
&\|(\nabla\times J_q)\otimes v_{q+1}+v_{q+1}\otimes (\nabla\times J_q)\|_{L^p}\\
\lesssim &\|\nabla\times J_q\|_{L^\infty}\|v_{q+1}\|_{L^p}\\
\lesssim &\ \lambda_q^3\left(\lambda_{q+1}^{-1}\delta_{q+1}^{\frac12}\ell^{-\frac12(1-\frac1p)}r^{\frac32-\frac3p}+\ell^{-1}\mu^{-1}(\lambda_{q+1}\sigma)^{-1}\delta_{q+1}r^{3-\frac3p}\right)\\
\lesssim &\ \lambda_q^3\ell^{-1}\mu^{-1}(\lambda_{q+1}\sigma)^{-1}\delta_{q+1}r^{3-\frac3p};
\end{split}
\end{equation}
and similarly, by (\ref{b-cn}) and (\ref{w-1p}),
\begin{equation}\notag
\begin{split}
&\|(\nabla\times w_{q+1})\otimes B_{q}+B_{q}\otimes (\nabla\times w_{q+1})\|_{L^p}\\
\lesssim &\|B_q\|_{L^\infty}\|w_{q+1}\|_{W^{1,p}}\\
\lesssim &\ \lambda_q^3\ell^{-2}\lambda_{q+1}r^{\frac32-\frac3p}.
\end{split}
\end{equation} 
Combining (\ref{Jn-p}), (\ref{v-pp})-(\ref{v-tp}), (\ref{v-w1p}), (\ref{vt-w1p}), (\ref{b-cn}), (\ref{w-p}), and (\ref{w-1p}) 
yields
\begin{equation}\notag
\begin{split}
&\|\nabla(J_q\times v_{q+1})+\nabla(w_{q+1}\times B_q)\|_{L^p}\\
\lesssim &\|\nabla J_q\|_{L^\infty}\|v_{q+1}\|_{L^p}+\|J_q\|_{L^\infty}\|v_{q+1}\|_{W^{1,p}}\\
&+\|\nabla B_q\|_{L^\infty}\|w_{q+1}\|_{L^p}+\|B_q\|_{L^\infty}\|w_{q+1}\|_{W^{1,p}}\\
\lesssim &\ \lambda_q^3 \ell^{-1}\mu^{-1}\lambda_{q+1}^{-1}\sigma^{-1}r_{q+1}^{3-\frac3p}
+\lambda_q^3\left(\ell^{-2}r^{\frac32-\frac3p}+\mu^{-1}\delta_{q+1}\ell^{-1}r^{4-\frac3p}\right)\\
&+\lambda_q^4\left(\delta_{q+1}^{\frac12}\ell^{-\frac12(1-\frac1p)}r^{\frac32-\frac3p}+\ell^{-2}\lambda_{q+1}r^{\frac32-\frac3p}\right).
\end{split}
\end{equation} 
Other terms in $R_{\mathrm{linear}}$ can be estimated similarly.
Summarizing the estimates above and taking into account the choice of parameters (\ref{pp1})-(\ref{pp2}) concludes the proof.

\cbdu

\medskip

\subsubsection{Correction terms.} 

\begin{Lemma}\label{le-corrector}
For $p>1$ close enough to 1, and a sufficiently small constant $\varepsilon_R>0$ depending on $p$, the following estimates hold:
\begin{equation}\notag
\|R_{\{\mathrm{cor},2\}}\|_{L^p}+\|R_{\{\mathrm{cor},3\}}\|_{L^p}+\|\tilde M_{q+1}\|_{L^p}\lesssim  \lambda_{q+1}^{-2\varepsilon_R}\delta_{q+1}.
\end{equation}
\end{Lemma}
\pf
Recall from (\ref{eq-r-corrector2})
\begin{equation}\notag
\begin{split}
\nabla\cdot R_{\mathrm{corrector}}
=&  \nabla\cdot (R_{\{\mathrm{cor},2\}}+R_{\{\mathrm{cor},3\}})+\nabla\times\nabla\times \tilde M_{q+1}
+\nabla\tilde p_{q+1,1}
\end{split}
\end{equation}
with
\[\tilde M_{q+1}=w_{q+1}\times v_{q+1}^t,\]
and $R_{\{\mathrm{cor},2\}}$, $R_{\{\mathrm{cor},3\}}$ defined in (\ref{eq-r-corr-2}), (\ref{eq-r-corr-3}) respectively.
We can estimate $\tilde M_{q+1}$ as, in view of (\ref{w-p}), (\ref{v-tp}) and (\ref{pp1})-(\ref{pp2})
\begin{equation}\notag
\begin{split}
\|\tilde M_{q+1}\|_{L^p}\leq &\|w_{q+1}\|_{L^{2p}}\|v_{q+1}^t\|_{L^{2p}}\\
\lesssim &\ \delta_{q+1}^{\frac12}\ell^{-\frac12(1-\frac1p)}r^{\frac32-\frac3{2p}}\ell^{-1}\mu^{-1}\lambda_{q+1}^{-1}\sigma^{-1}\delta_{q+1}r^{3-\frac3{2p}}\\
\lesssim &\ \ell^{-3}\mu^{-1}\lambda_{q+1}^{-1}\sigma^{-1}\delta_{q+1}r^{\frac92-\frac3{p}}\\
\lesssim &\ \lambda_{q+1}^{-2\varepsilon_R}\delta_{q+1}.
\end{split}
\end{equation}
We turn to the estimates of $R_{\{\mathrm{cor},2\}}$, $R_{\{\mathrm{cor},3\}}$, and $\tilde M_{q+1,2}$ which are trivial. Following from (\ref{v-cp}) and (\ref{w-1p}), it has
\begin{equation}\notag
\begin{split}
\|R_{\{\mathrm{cor},2\}}\|_{L^p}\leq&\ \|w_{q+1}\|_{W^{1,2p}}\|v_{q+1}^c\|_{L^{2p}}\\
\lesssim&\ \lambda_{q+1}^{-1}\delta_{q+1}^{\frac12}\ell^{-\frac12(1-\frac1p)}\sigma r^{\frac52-\frac3{2p}}\ell^{-2}\lambda_{q+1}r^{\frac32-\frac3{2p}}\\
\lesssim&\ \ell^{-3}\delta_{q+1}^{\frac12}\sigma r^{4-\frac3{p}}\\
\lesssim &\ \lambda_{q+1}^{-2\varepsilon_R}\delta_{q+1}.
\end{split}
\end{equation}
By (\ref{w-1p}) and (\ref{v-pp}), we have, for $p>1$ sufficiently close to 1
\begin{equation}\notag
\begin{split}
\|R_{\{\mathrm{cor},3\}}\|_{L^p}\leq&\ \|w_{q+1}^c+w_{q+1}^t\|_{W^{1,2p}}\|v_{q+1}^p\|_{L^{2p}}\\
\lesssim&\ \lambda_{q+1}^{-1}\delta_{q+1}^{\frac12}\ell^{-\frac12(1-\frac1p)} r^{\frac32-\frac3{2p}}\ell^{-2}\lambda_{q+1}r^{\frac32-\frac3{2p}}\\
\lesssim&\ \ell^{-3}\delta_{q+1}^{\frac12} r^{3-\frac3{p}}\\
\lesssim &\ \lambda_{q+1}^{-2\varepsilon_R}\delta_{q+1},
\end{split}
\end{equation}
where the last step follows from the choice of the parameters (\ref{pp1})-(\ref{pp2}).

\cbdu

\medskip

\subsubsection{Oscillation terms.} 

\begin{Lemma}\label{le-oscill}
Recall from (\ref{eq-r-osc-rewrite})
\begin{equation}\notag
\nabla\cdot R_{\mathrm{oscillation}}= \nabla\cdot(R_{\{\mathrm{osc},1\}}+R_{\{\mathrm{osc},2\}})
\end{equation}
with $R_{\{\mathrm{osc},1\}}$ and $R_{\{\mathrm{osc},2\}}$ respectively defined in (\ref{eq-r-osc-1}) and (\ref{eq-r-osc-2}). For $p>1$ sufficiently close to 1 and an arbitrarily small constant $\varepsilon_R>0$, we have
\begin{equation}\notag
\|R_{\{\mathrm{osc},1\}}\|_{L^p}+\|R_{\{\mathrm{osc},2\}}\|_{L^p}\lesssim \lambda_{q+1}^{-2\varepsilon_R}\delta_{q+1}.
\end{equation}
\end{Lemma}
\pf
Notice that $\nabla\cdot R_{\{\mathrm{osc},1\}}$ is a similar oscillation term as $\nabla\cdot \tilde R_{\mathrm{oscillation}}$ for the NSE in \cite{BV}, and hence can be estimated in an analogous way.  Without showing details, we claim
\begin{equation}\label{est-r-oscill-first3}
\|R_{\{\mathrm{osc},1\}}\|_{L^p}\lesssim \lambda_{q+1}^{-2\varepsilon_R}\delta_{q+1}.
\end{equation}

We estimate $R_{\{\mathrm{osc},2\}}$ as follows. 
Applying (\ref{w-p}), (\ref{est-We-sobolev}), and  (\ref{pp1})-(\ref{pp2}), we deduce that for $p>1$ close enough to $1$
\begin{equation}\notag
\begin{split}
&\|w_{q+1}^p\otimes \mathbb W_{\epsilon,1}+\mathbb W_{\epsilon,1}\otimes w_{q+1}^p\|_{L^p}\\
\lesssim & \|w_{q+1}\|_{L^p}\|\mathbb W_{\epsilon,1}\|_{L^\infty}\\
\lesssim &\delta_{q+1}^{\frac12}\ell^{-\frac12(1-\frac1p)} r^{\frac32-\frac3p}\delta_{q+1}^{\frac12}\ell^{-\frac12}\sigma r^{\frac52}\\
\lesssim&\lambda_{q+1}^{-2\varepsilon_R}\delta_{q+2}.
\end{split}
\end{equation}
Using (\ref{w-p}), (\ref{v-pp}) and (\ref{pp1})-(\ref{pp2}) gives us 
\begin{equation}\notag
\begin{split}
&\|w_{q+1}^c\otimes v_{q+1}^p+v_{q+1}^p\otimes w_{q+1}^c\|_{L^p}\\
\lesssim & \|w_{q+1}^c\|_{L^{2p}} \|v_{q+1}^p\|_{L^{2p}}\\
\lesssim & \delta_{q+1}^{\frac12}\ell^{-\frac12(1-\frac1{2p})}r^{\frac32-\frac3{2p}} \lambda_{q+1}^{-1}\delta_{q+1}^{\frac12}\ell^{-\frac12(1-\frac1{2p})} r^{\frac32-\frac3{2p}}\\
\lesssim&\lambda_{q+1}^{-2\varepsilon_R}\delta_{q+2}
\end{split}
\end{equation}
for $p>1$ sufficiently close to $1$. Applying (\ref{est-We-sobolev}) leads to
\begin{equation}\notag
\|\mathbb W_{\epsilon,1}\otimes \mathbb W_{\epsilon,1}\|_{L^p}
\leq \|\mathbb W_{\epsilon,1}\|_{L^{2p}}^2
\lesssim \left(\delta_{q+1}^{\frac12}\ell^{-\frac12(1-\frac1{2p})}\sigma r^{\frac52-\frac3{2p}}\right)^2
\lesssim\lambda_{q+1}^{-2\varepsilon_R}\delta_{q+2}
\end{equation}
where the last step holds for $p>1$ close enough to $1$ in view of (\ref{pp1})-(\ref{pp2}). 
In the end, employing (\ref{est-We-sobolev}), (\ref{v-pp}) and (\ref{pp1})-(\ref{pp2}) leads to 
\begin{equation}\notag
\begin{split}
&\|\mathcal R\nabla\times(\mathbb W_{\epsilon,1}\nabla v_{q+1}^p\|_{L^p}\\
\lesssim &\|\mathbb W_{\epsilon,1}\|_{L^\infty} \| \nabla v_{q+1}^p\|_{L^p}\\
\lesssim & \delta_{q+1}^{\frac12}\ell^{-\frac12}\sigma r^{\frac52} \delta_{q+1}^{\frac12}\ell^{-\frac12(1-\frac1{p})} r^{\frac32-\frac3{p}}\\
\lesssim&\lambda_{q+1}^{-2\varepsilon_R}\delta_{q+2}
\end{split}
\end{equation}
for $p>1$ close enough to $1$. 
The last inequalities together with equation (\ref{eq-r-osc-2}) yield 
\begin{equation}\label{est-r-oscill-middle}
\|R_{\{\mathrm{osc},2\}}\|_{L^p}\lesssim \lambda_{q+1}^{-2\varepsilon_R}\delta_{q+1}.
\end{equation}
The conclusion of the lemma follows immediately from (\ref{est-r-oscill-first3}) and (\ref{est-r-oscill-middle}). 

\cbdu

\medskip

\subsection{The energy iteration}
\label{sec-energy}

\begin{Lemma}\label{le-energy3}
If $\rho_0(t)\neq 0$, then the energy of the current density $J_{q+1}$ satisfies 
\begin{equation}\notag
\left| E(t)-\int_{\mathbb T^3}|J_{q+1}(x,t)|^2\, \mathrm dx-\frac{\delta_{q+2}}{2}\right|\leq \frac{\delta_{q+2}}{4}.
\end{equation}
\end{Lemma}

\begin{Lemma}\label{le-energy4}
If $\rho_0(t)=0$, then $J_{q+1}(\cdot, t)\equiv 0$, $R_{q+1}(\cdot, t)\equiv 0$ and
\begin{equation}\notag
E(t)-\int_{\mathbb T^3}|J_{q+1}(x,t)|^2\, \mathrm dx\leq \frac34\delta_{q+2}.
\end{equation}
\end{Lemma}
The proof of Lemma \ref{le-energy3} and Lemma \ref{le-energy4} follows closely as the proof of Lemma 6.2 and Lemma 6.3 in \cite{BV}. The two estimates in Lemma \ref{le-energy3} and Lemma \ref{le-energy4} immediately implies (\ref{energy-q}) for $q+1$. On the other hand, if 
\begin{equation}\notag
E(t)-\int_{\mathbb T^3}|J_{q+1}(x,t)|^2\,\mathrm dx\leq \frac{\delta_{q+2}}{100},
\end{equation}
it follows from Lemma \ref{le-energy3} that $\rho_0(t)=0$. Thus, Lemma \ref{le-energy4} guarantees $J_{q+1}(t)\equiv 0$ and $R_{q+1}(t)\equiv 0$, which shows (\ref{energy-q1}) for $q+1$.

Now we can conclude that the proof of Proposition \ref{le-iterative} is complete.

\bigskip

\section{Non-uniqueness of the Hall MHD system}
\label{sec-hmhd}

In this section, we come back to the 3D Hall-MHD system (\ref{HMHD}) with $\zeta=1$ and demonstrate that non-unique Leray-Hopf weak solutions can be actually constructed for this coupled system of the NSE and the Hall equation. That is, we prove Theorem \ref{thm}.

We consider the approximating system 
\begin{equation}\label{hmhd-q}
\begin{split}
\partial_t u_q+(u_q\cdot\nabla) u_q+\nabla p_q=&\ \Delta u_q+(B_q\cdot\nabla) B_q,\\
\partial_t J_q+\nabla\times\nabla\times (B_{q}\times u_{q}) +& \nabla\times\nabla\times (J_q\times B_q)\\
=&\ \Delta J_q+\nabla\cdot R_q^s,\\
\nabla\cdot u_q=&\ 0.
\end{split}
\end{equation}
The plan is to apply convex integration framework only to the equation of the current density $J_q$ and solve the NSE with force $(B_q\cdot\nabla) B_q$ at every level of the convex integration. The detailed scheme is described below:
\begin{itemize}
\item [(i)]
Start with the trivial choice of $(u_0, B_0, J_0, R_0^s)=(0, 0, 0, 0)$ which satisfies (\ref{hmhd-q}) automatically.
\item [(ii)]
Construct appropriate (small enough) perturbations $v_1^s=B_1-B_0$ and $w_1^s=J_1-J_0$ to obtain $J_1=J_0+w^s_1$ and $B_1=B_0+v^s_1$.
\item [(iii)]
Solve the NSE with force $(B_1\cdot\nabla) B_1$,
\begin{equation}\label{eq-nse1}
\partial_t u_1+(u_1\cdot\nabla)u_1+\nabla p_1=\Delta u_1+(B_1\cdot\nabla) B_1
\end{equation}
and obtain $u_1$.
\item [(iv)]
With $B_1, J_1$, and $u_1$, we derive the new stress tensor $R_1^s$ which satisfies 
\[\partial_t J_1+\nabla\times \nabla \times (B_1\times u_1)+\nabla\times \nabla \times (J_1\times B_1)=\Delta J_1+\nabla\cdot R^s_{1}.\]
We have, equivalently, by noticing $u_0= B_0= J_0= R_0^s=0$, $B_1=v_1^s$ and  $J_1=w_1^s$, 
\[
\begin{split}
\nabla\cdot R^s_{1}= &\ \partial_t w_{1}-\Delta w_{1}+\nabla\times \nabla\times (J_{0}\times v_{1})\\
&+\nabla\times \nabla\times (w_{1}\times B_{0})+\nabla\times \nabla\times (w_{1}\times v_{1})\\
&+\nabla\cdot R_{0}^s+ \nabla\times\nabla\times (B_{1}\times u_{1}-B_0\times u_0).
\end{split}
\]

For such stress tensor $R^s_{1}$, we know the tuplet $(u_1, B_1, J_1, R_1^s)$ satisfies system (\ref{hmhd-q}) at level $q=1$. It completes the first iteration.
\item [(v)]
To correct the stress tensor $R_1^s$,  we move forward to the second iteration which also contains a few steps. First, we construct increments $v^s_2=B_2-B_1$ and $w^s_2=J_2-J_1$ and hence obtain $B_2=B_1+v^s_2$ and $J_2=J_1+w^s_2$.
\item [(vi)]
We then solve the NSE with force $(B_2\cdot\nabla) B_2$,
\[\partial_t u_2+(u_2\cdot\nabla)u_2+\nabla p_2=\Delta u_2+(B_2\cdot\nabla) B_2\]
to obtain $u_2$. 
\item [(vii)]
With $B_2, J_2$, and $u_2$, we derive the stress tensor $R_2^s$ satisfying 
\[\partial_t J_2+\nabla\times \nabla \times (B_2\times u_2)+\nabla\times \nabla \times (J_2\times B_2)=\Delta J_2+\nabla\cdot R^s_2.\]
Thus, we have
\begin{equation}\notag
\begin{split}
\nabla\cdot R^s_2=&\ \partial_t J_2-\Delta J_2+\nabla\times \nabla \times (B_2\times u_2)+\nabla\times \nabla \times (J_2\times B_2)\\
=&\ \partial_t (J_1+w^s_2)-\Delta (J_1+w^s_2)+\nabla\times \nabla \times (B_2\times u_2)\\
&+\nabla\times \nabla \times ((J_1+w^s_2)\times (B_1+v^s_2)).
\end{split}
\end{equation}
It follows from the equation of $J_1$ in Step (iv) that
\[\partial_t J_1- \Delta J_1+\nabla\times \nabla \times (J_1\times B_1)=\nabla\cdot R^s_{1}-\nabla\times \nabla \times (B_1\times u_1).\]
Combining the last two equations leads to
\begin{equation}\notag
\begin{split}
\nabla\cdot R^s_2=&\ \partial_t w^s_2-\Delta w^s_2+\nabla\times \nabla\times (J_{1}\times v^s_{2})\\
&+\nabla\times \nabla\times (w^s_{2}\times B_{1})+\nabla\times \nabla\times (w^s_{2}\times v^s_{2})\\
&+\nabla\cdot R_{1}^s+ \nabla\times\nabla\times (B_{2}\times u_{2}-B_1\times u_1).
\end{split}
\end{equation}
Thus, we complete the second iteration and obtain $(u_2, B_2, J_2, R_2^s)$ with $R_2^s$ satisfying the last equation.
\item [(viii)]
Repeat Steps (v)-(vii) iteratively to obtain a sequence $\{(u_{q+1}, B_{q+1}, J_{q+1}$, $ R_{q+1}^s)\}$ satisfying (\ref{hmhd-q}) and the stress tensor $R_{q+1}^s$ satisfies (\ref{eq-hmhd-r-q}) in Lemma \ref{le-hmhd-new-stress} below.
\item [(viiii)]
Prove that the sequence $\{(u_{q+1}, B_{q+1}, J_{q+1}, R_{q+1}^s)\}$ converges to $(u, B, J, 0)$ with functions $u, B, J$ satisfying
\begin{equation}\notag
J=\nabla\times B, \ \ u\in L^\infty(L^2)\cap L^2(H^1), \ \ B\in L^\infty(L^2)\cap L^2(H^1),
\end{equation}
and $(u,B)$ is a weak solution of the Hall-MHD system (\ref{HMHD}).
\end{itemize}

\begin{Remark}\label{rk-comparison1}
We notice that at level $q=1$, the stress tensor $R_1^s$ for the system (\ref{hmhd-q}) and the stress tensor $R_1$ for the equation (\ref{eq-div-q-alt}) are different. Indeed, we have
\begin{equation}\notag
\begin{split}
\nabla\cdot R_1^s=&\ \partial_t J_1-\Delta J_1+\nabla\times \nabla \times (B_1\times u_1)+\nabla\times \nabla \times (J_1\times B_1),\\
\nabla\cdot R_1=&\ \partial_t J_1-\Delta J_1+\nabla\times \nabla \times (J_1\times B_1).
\end{split}
\end{equation}
In principle, the increments $(v^s_2, w^s_2)$ (or $(v_2, w_2)$) are constructed to correct the stress tensor $R_1^s$ (or $R_1$). Thus, due to the difference of $R_1^s$ and $R_1$, we should have $(v^s_2, w^s_2)\neq (v_2, w_2)$. However, $(v^s_2, w^s_2)$ can be constructed analogously as $(v_2, w_2)$ by using intermittent Beltrami flows. The difference relies on the coefficients $a_{\xi}$ of the Beltrami flows, which depend on the stress tensor, see details of Subsections \ref{sec-design} and \ref{sec-est-blocks}. 
\end{Remark}

\begin{Remark}\label{rk-comparison2}
The increments $(v^s_{q+1}, w^s_{q+1})$ for $q\geq0$ can be constructed similarly as $(v_{q+1}, w_{q+1})$ such that the estimates in Lemma \ref{le-vq} and Lemma \ref{le-w} hold for $(v^s_{q+1}, w^s_{q+1})$ as well. Note that the coefficients $a_{\xi}$ of the Beltrami flows in the construction of $(v^s_{q+1}, w^s_{q+1})$ depend on the stress tensor $R^s_q$.
\end{Remark}

\begin{Remark}\label{rk-comparison3}
In the iterating process, the stress tensor $R^s_{q+1}$ for $q\geq0$ satisfies equation (\ref{eq-hmhd-r-q}) below, i.e.
\[\nabla\cdot R_{q+1}^s= \nabla\cdot \tilde R_{q+1}^s+\nabla\times\nabla\times M^\epsilon_{q+1}\]
with
\begin{equation}\notag
\begin{split}
\nabla\cdot \tilde R_{q+1}^s=&\ \nabla\cdot R_{q}^s+\partial_t w^s_{q+1}-\Delta w^s_{q+1}+\nabla\times \nabla\times (J_{q}\times v^s_{q+1})\\
&+\nabla\times \nabla\times (w^s_{q+1}\times B_{q})+\nabla\times \nabla\times (w^s_{q+1}\times v^s_{q+1}),\\
M^\epsilon_{q+1}= &\ B_{q+1}\times u_{q+1}-B_q\times u_q.
\end{split}
\end{equation} 
We observe that $\nabla\cdot \tilde R_{q+1}^s$ has the same structure as $\nabla\cdot R_{q+1}$ of (\ref{eq-emhd-R-q4}). Hence, $\tilde R_{q+1}^s$ can be dealt with and estimated in an analogous way as that of $R_{q+1}$. We would like to point out how to deal with $\nabla\times\nabla\times M^\epsilon_{q+1}$, for instance in the first iteration as $q=0$: 
\begin{equation}\notag
\begin{split}
\nabla\cdot R_{1}^s=&\ \nabla\cdot \tilde R_{1}^s+\nabla\times\nabla\times M^\epsilon_{1}, \\
M^\epsilon_{1}= &\ B_{1}\times u_{1}-B_0\times u_0= B_{1}\times u_{1}
\end{split}
\end{equation}
which is also reflected in Remark \ref{rk-comparison1}. We notice that it is sufficient to estimate $M^\epsilon_{1}$ rather than $\div^{-1} \nabla\times\nabla\times M^\epsilon_{1}$; indeed, in the process of passing limit in 
\[\int_{\mathbb T^3} (\nabla\cdot R_{q}^s) \cdot \varphi \, dx= \int_{\mathbb T^3} (\nabla\cdot \tilde R_{q}^s+\nabla\times\nabla\times M^\epsilon_{q}) \cdot \varphi \, dx\]
as $q\to \infty$, we can use integration by parts to move derivatives to the test function $\varphi$. We construct $v_1^s=B_1$ such that $\|B_1\|_{L^\infty (0,T; L^2(\mathbb T^3))}\leq c\delta_2$ for a small constant $c$ and $B_1\otimes B_1\in L^2(0,T; L^2(\mathbb T^3))$. Thus, the force in the NSE (\ref{eq-nse1}) satisfies 
\[(B_1\cdot \nabla) B_1=\nabla\cdot (B_1\otimes B_1)\in L^2(0,T; W^{-1, 2}(\mathbb T^3));\]
consequently, the solution $u_1$ of (\ref{eq-nse1}) belongs to the space $L^\infty(0,T; L^2(\mathbb T^3))$, see \cite{Tem}. Therefore, we have
\begin{equation}\notag
\begin{split}
\| M^\epsilon_{1}\|_{L^\infty (0,T; L^1(\mathbb T^3))}=&\ \|B_{1}\times u_{1}\|_{L^\infty (0,T; L^1(\mathbb T^3))}\\
\leq&\ \|B_{1}\|_{L^\infty (0,T; L^2(\mathbb T^3))} \|u_{1}\|_{L^\infty (0,T; L^2(\mathbb T^3))}\\
\leq&\ c \delta_2.
\end{split}
\end{equation}
\end{Remark}

The scheme described in (i)-(viiii) involves two major bulks:  solving the NSE of $u_q$ and applying convex integration on the $J_q$ equation. Details are demonstrated by proving the following iterative argument.

\begin{Proposition}\label{hmhd-iterative}
There exists an absolute constant $C>0$ and a sufficiently small parameter $\varepsilon_R$ depending on $b$ and $\beta$ such that the following inductive statement holds.
Let $(u_q, p_q, B_q, J_q, R_q^s)$ be a solution of the approximating equation (\ref{hmhd-q}) on $\mathbb T^3\times [0,T]$ satisfying:
\begin{equation}\label{bq-c1-h}
\|B_q\|_{C^1_{x,t}}\leq \lambda_q^3,
\end{equation}
\begin{equation} \label{jq-c1-h}
\|J_q\|_{C^1_{x,t}}\leq \lambda_q^4, 
\end{equation}
\begin{equation}\label{energy-q-h}
0\leq E(t)-\int_{\mathbb T^3}|J_q|^2\, \mathrm dx\leq \delta_{q+1},
\end{equation}
and 
\begin{equation}\label{energy-q1-h}
E(t)-\int_{\mathbb T^3}|J_q|^2\, \mathrm dx\leq \frac{\delta_{q+1}}{100} \ \ \mbox {implies} \ \ J_{q}(\cdot, t)\equiv 0 \ \mbox{and} \ R_q^s(\cdot, t)\equiv 0.
\end{equation}
In addition, we assume 
\begin{equation}\label{RM-h}
\nabla\cdot R_q^s=\nabla\cdot \tilde R_q^s+\nabla\times\nabla\times  M_q^\epsilon
\end{equation}
with $\tilde R_q^s$ being a symmetric traceless stress tensor 
and $ M_q^\epsilon$ being a vector field
which satisfy
\begin{equation}\label{est-rqs-tilde}
\|\tilde R_q^s\|_{L^\infty(L^1)}\leq \lambda_q^{-\varepsilon_R}\delta_{q+1}, 
\end{equation}
\begin{equation}\label{est-increm-M}
\|M_q^\epsilon\|_{L^\infty(L^1)}\leq \lambda_q^{-\varepsilon_R}\delta_{q+1}+C\|z_{q}\|_{L^2}.
\end{equation}

Then we can find another solution $(u_{q+1}, p_{q+1}, B_{q+1}, J_{q+1}, R_{q+1}^s)$ of (\ref{hmhd-q}) satisfying (\ref{bq-c1-h})-(\ref{est-increm-M}) with $q$ replaced by $q+1$. Moreover, 
the increments $v^s_{q+1}=B_{q+1}-B_q$, $w^s_{q+1}=J_{q+1}-J_q$ and $z^s_{q+1}=u_{q+1}-u_q$ satisfy 
\begin{equation}\label{est-increm-h}
\|v^s_{q+1}\|_{L^2}\leq C\lambda_{q+1}^{-1}\delta_{q+1}^{1/2}, \ \ \|w^s_{q+1}\|_{L^2}\leq C\delta_{q+1}^{1/2},
\end{equation}
\begin{equation}\label{est-increm-u}
\lim_{q\to\infty}\|z^s_{q+1}\|_{L^p}=0, \ \ 1\leq p\leq 2.
\end{equation}
\end{Proposition}

In analogy with Proposition \ref{le-iterative} and Theorem \ref{thm-h}, a proof of Theorem \ref{thm} follows immediately from Proposition \ref{hmhd-iterative}; thus the details are omitted. 

In order to prove Proposition \ref{hmhd-iterative}, we adapt the same construction of perturbations $v^s_{q+1}=B_{q+1}-B_q$ and $w^s_{q+1}=J_{q+1}-J_q$ as that of $v_{q+1}$ and $w_{q+1}$ for the Hall equation in Section \ref{sec-hall}, however, with different coefficients for the Beltrami flows, as indicated in Remark \ref{rk-comparison2}. Since $v^s_{q+1}$ and $w^s_{q+1}$ also satisfy the estimates in Lemma \ref{le-vq} and Lemma \ref{le-w}, respectively, we claim that (\ref{bq-c1-h}), (\ref{jq-c1-h}), and (\ref{est-increm-h}) hold for $q+1$.

We continue to complete the proof of Proposition \ref{hmhd-iterative} in Subsections \ref{sec-nse-q}, \ref{sec-hmhd-new-stress} and \ref{sec-other-estimate} below.

\medskip

\subsection{Weak solution $u_{q+1}$ of the NSE in $L^\infty(L^2)\cap L^2(H^1)$}
\label{sec-nse-q}
We consider the forced NSE 
\begin{equation}\label{nse-q}
\partial_tu_{q+1}+(u_{q+1}\cdot\nabla) u_{q+1}+\nabla p_{q+1}=\Delta u_{q+1}+\nabla\cdot(B_{q+1}\otimes B_{q+1}).
\end{equation}
By construction, we have
\begin{equation}\notag
B_{q+1}=B_0+\sum_{j=0}^{j=q}v^s_{j+1}, \ \ J_{q+1}=J_0+\sum_{j=0}^{j=q}w^s_{j+1}
\end{equation}
with $\|v^s_{q+1}\|_{L^2}\leq C\lambda_{q+1}^{-1}\delta_{q+1}^{1/2}$ and $\|w^s_{q+1}\|_{L^2}\leq C\delta_{q+1}^{1/2}$. It is then obvious that $\|B_{q+1}\|_{L^2}\leq C$ and $\|J_{q+1}\|_{L^2}\leq C$ which implies 
$B_{q+1}\in L^\infty(0,T; H^1(\mathbb T^3))$, 
since $B_{q+1}$ is divergence free. 

It follows from the Sobolev embedding theorem that $B_{q+1}\otimes B_{q+1}$ is in the space
$L^2(0,T; L^3(\mathbb T^3))$, and hence in $L^2(0,T; L^2(\mathbb T^3))$ as well. Thus we have $\nabla\cdot (B_{q+1}\otimes B_{q+1})\in L^2(0,T; W^{-1, 2})$. Then there exists a weak solution $u_{q+1}$ of (\ref{nse-q}) with $u_{q+1}\in L^\infty(0,T; L^2(\mathbb T^3))\cap L^2(0,T; H^1(\mathbb T^3))$, see \cite{Tem}. 

Upon writing $u_{q+1}$ as the sum of increments,
\begin{equation}\notag
u_{q+1}=u_0+\sum_{j=0}^{j=q}(u_{j+1}-u_j)=u_0+\sum_{j=0}^{j=q}z^s_{j+1},
\end{equation}
the fact $u_{q+1}\in L^\infty(0,T; L^2(\mathbb T^3))$ implies 
\begin{equation}\notag
\lim_{q\to \infty}\|z^s_{q+1}(t)\|_{L^2(\mathbb T^3)}=0, \ t\in[0,T].
\end{equation}
Moreover, we have $\lim_{q\to \infty}\|z^s_{q+1}(t)\|_{L^p(\mathbb T^3)}=0$, $0\leq t\leq T$, for all $p\in[1,2]$. Therefore, (\ref{est-increm-u}) is justified.

\medskip

\subsection{The new stress tensor $R_{q+1}^s$ and its estimate}
\label{sec-hmhd-new-stress}

With $v^s_{q+1}=B_{q+1}-B_q$ and $w^s_{q+1}=J_{q+1}-J_q$ constructed following the same line of Section \ref{sec-hall} and $u_{q+1}$ obtained in Section \ref{sec-nse-q}, we proceed to derive the stress tensor $R_{q+1}^s$ satisfying the equation of $J_{q+1}$ in (\ref{hmhd-q}) at the level $q+1$. 
Compared to the $J_q$ equation in (\ref{eq-div-q}), there is one extra term 
 $\nabla\times\nabla\times (B_q\times u_q)$ in the $J_q$ equation of (\ref{hmhd-q}). Thus, $R_{q+1}^s$ will be different from $R_{q+1}$ mainly due to the interaction of this extra nonlinear term. In the following, we will show that: 
 \begin{Lemma}\label{le-hmhd-new-stress}
 Let $v^s_{q+1}=B_{q+1}-B_q$ and $w^s_{q+1}=J_{q+1}-J_q$ be the increments appropriately constructed similarly as in Section \ref{sec-hall}, which depend on $R^s_q$. Let $u_{q+1}$ be the new velocity obtained in Subsection \ref{sec-nse-q}. We choose the symmetric traceless stress tensor $R_{q+1}^s$ such that
\begin{equation}\label{eq-hmhd-r-q}
\begin{split}
\nabla\cdot R_{q+1}^s=&\ \partial_t w^s_{q+1}-\Delta w^s_{q+1}+\nabla\times \nabla\times (J_{q}\times v^s_{q+1})\\
&+\nabla\times \nabla\times (w^s_{q+1}\times B_{q})+\nabla\times \nabla\times (w^s_{q+1}\times v^s_{q+1})\\
&+\nabla\cdot R_{q}^s+ \nabla\times\nabla\times (B_{q+1}\times u_{q+1}-B_q\times u_q)\\
=:&\ \nabla\cdot \tilde R_{q+1}^s+\nabla\times\nabla\times M^\epsilon_{q+1}
\end{split}
\end{equation} 
with 
\[M^\epsilon_{q+1}= B_{q+1}\times u_{q+1}-B_q\times u_q\] 
and $\nabla\cdot \tilde R_{q+1}^s$ being the rest terms of $\nabla\cdot R_{q+1}^s$. 
Then, the tuplet $(u_{q+1}, B_{q+1}$, $ J_{q+1}, R_{q+1}^s)$ solves system (\ref{hmhd-q}) with $q$ replaced by $q+1$.
\end{Lemma}
\pf
Since
\[B_{q+1}=B_q+v^s_{q+1}, \ \ \ J_{q+1}=J_q+w^s_{q+1},\]
we have
\[\partial_t J_{q+1}=\partial_t J_q+\partial_t w^s_{q+1}, \ \ \ \Delta J_{q+1}=\Delta J_q+\Delta w^s_{q+1},\]
and 
\begin{equation}\notag
\begin{split}
&\nabla\times \nabla\times (J_{q+1}\times B_{q+1})\\
=& \nabla\times \nabla\times ((J_{q}+w^s_{q+1})\times (B_{q}+v^s_{q+1}))\\
=& \nabla\times \nabla\times (J_{q}\times B_{q})+\nabla\times \nabla\times (J_{q}\times v^s_{q+1})\\
&+\nabla\times \nabla\times (w^s_{q+1}\times B_{q})+\nabla\times \nabla\times (w^s_{q+1}\times v^s_{q+1}).
\end{split}
\end{equation}
Therefore, we continue to deduce
\begin{equation}\label{eq-j-q1}
\begin{split}
&\partial_t J_{q+1}+ \nabla\times \nabla\times (B_{q+1}\times u_{q+1})+\nabla\times \nabla\times (J_{q+1}\times B_{q+1})\\
=& \partial_t J_q+\partial_t w^s_{q+1}+ \nabla\times \nabla\times (B_{q+1}\times u_{q+1})\\
&+ \nabla\times \nabla\times (J_{q}\times B_{q})+\nabla\times \nabla\times (J_{q}\times v^s_{q+1})\\
&+\nabla\times \nabla\times (w^s_{q+1}\times B_{q})+\nabla\times \nabla\times (w^s_{q+1}\times v^s_{q+1}).
\end{split}
\end{equation}
Recall that the tuplet $(u_q, B_q, J_q, R_q^s)$ satisfies equation (\ref{hmhd-q}). Hence, 
\begin{equation}\label{eq-j-q0}
\begin{split}
&\partial_t J_q+ \nabla\times \nabla\times (J_{q}\times B_{q})\\
=&\Delta J_q-\nabla\times \nabla\times (B_{q}\times u_{q})+\nabla\cdot R_q^s.
\end{split}
\end{equation}
Inserting (\ref{eq-j-q0}) into the right hand side of (\ref{eq-j-q1}) yields
\begin{equation}\label{eq-j-q2}
\begin{split}
&\partial_t J_{q+1}+ \nabla\times \nabla\times (B_{q+1}\times u_{q+1})+\nabla\times \nabla\times (J_{q+1}\times B_{q+1})\\
=& \Delta J_q-\nabla\times \nabla\times (B_{q}\times u_{q})+\nabla\cdot R_q^s +\partial_t w^s_{q+1}\\
&+ \nabla\times \nabla\times (B_{q+1}\times u_{q+1})
+\nabla\times \nabla\times (J_{q}\times v^s_{q+1})\\
&+\nabla\times \nabla\times (w^s_{q+1}\times B_{q})+\nabla\times \nabla\times (w^s_{q+1}\times v^s_{q+1}).
\end{split}
\end{equation}
Substituting $\Delta J_q$ by $\Delta J_{q+1}-\Delta w^s_{q+1}$ in (\ref{eq-j-q2}) shows that
\begin{equation}\label{eq-j-q3}
\begin{split}
&\partial_t J_{q+1}+ \nabla\times \nabla\times (B_{q+1}\times u_{q+1})+\nabla\times \nabla\times (J_{q+1}\times B_{q+1})\\
=&\Delta J_{q+1}-\Delta w^s_{q+1}-\nabla\times \nabla\times (B_{q}\times u_{q})+\nabla\cdot R_q^s +\partial_t w^s_{q+1}\\
&+ \nabla\times \nabla\times (B_{q+1}\times u_{q+1})
+\nabla\times \nabla\times (J_{q}\times v^s_{q+1})\\
&+\nabla\times \nabla\times (w^s_{q+1}\times B_{q})+\nabla\times \nabla\times (w^s_{q+1}\times v^s_{q+1}).
\end{split}
\end{equation}
Thus we choose $R_{q+1}^s$ such that 
\begin{equation}\notag
\begin{split}
\nabla\cdot R_{q+1}^s=&\ \partial_t w^s_{q+1}-\Delta w^s_{q+1}+\nabla\times \nabla\times (J_{q}\times v^s_{q+1})\\
&+\nabla\times \nabla\times (w^s_{q+1}\times B_{q})+\nabla\times \nabla\times (w^s_{q+1}\times v^s_{q+1})\\
&+\nabla\cdot R_{q}^s+ \nabla\times\nabla\times (B_{q+1}\times u_{q+1}-B_q\times u_q),
\end{split}
\end{equation}
We denote $M_{q+1}^\epsilon=B_{q+1}\times u_{q+1}-B_q\times u_q$, and choose $\tilde R_{q+1}^s$ such that 
\begin{equation}\label{eq-tilde-rq1}
\begin{split}
\nabla\cdot \tilde R_{q+1}^s=&\ \nabla\cdot R_{q}^s+\partial_t w^s_{q+1}-\Delta w^s_{q+1}+\nabla\times \nabla\times (J_{q}\times v^s_{q+1})\\
&+\nabla\times \nabla\times (w^s_{q+1}\times B_{q})+\nabla\times \nabla\times (w^s_{q+1}\times v^s_{q+1}).\\
\end{split}
\end{equation}
Hence the tuplet $(u_{q+1}, B_{q+1}$, $ J_{q+1}, R_{q+1}^s)$ satisfies (\ref{hmhd-q}) at the level $q+1$ thanks to (\ref{eq-j-q3}).

\cbdu

In the following, we estimate $\tilde R_{q+1}^s$, $M_{q+1}^\epsilon$ and hence $R_{q+1}^s$.
\begin{Lemma} \label{le-est-tilde-rq1s}
For $p>1$ sufficiently close to $1$ and $\varepsilon_R>0$ sufficiently small, the stress tensor $\tilde R_{q+1}^s$ satisfies 
\begin{equation}\notag
\|\tilde R_{q+1}^s\|_{L^p}\leq \lambda_{q+1}^{-\varepsilon_R}\delta_{q+2}.
\end{equation}
\end{Lemma}

\pf
As pointed out in Remark \ref{rk-comparison3}, $\tilde R_{q+1}^s$ of (\ref{eq-tilde-rq1}) has the same structure as the stress tensor $R_{q+1}$. In analogy of handling $R_{q+1}$ of (\ref{eq-emhd-R-q4}) in Subsection \ref{sec-emhd-rq1},  we can rewrite $\nabla\cdot \tilde R_{q+1}^s$ into
\begin{equation}\notag
\begin{split}
\nabla\cdot \tilde R_{q+1}^s=\nabla\cdot R^s_{\mathrm{linear}}+\nabla\cdot R^s_{\mathrm{corrector}}+\nabla\cdot R^s_{\mathrm{oscillation}}
\end{split}
\end{equation}
with
\begin{equation}\notag
\begin{split}
R^s_{\mathrm{linear}}=&\ R_{\mathrm{linear}}(v^s_{q+1}, w^s_{q+1}), \\ 
R^s_{\mathrm{corrector}}=&\ R_{\mathrm{corrector}}(v^s_{q+1}, w^s_{q+1}), \\
R^s_{\mathrm{oscillation}}=&\ R_{\mathrm{oscillation}}(v^s_{q+1}, w^s_{q+1}, R^s_q),
\end{split}
\end{equation}
where $R_{\mathrm{linear}}(v^s_{q+1}, w^s_{q+1})$ denotes the linear part of the stress tensor as in (\ref{eq-emhd-R-q6}) with $(v_{q+1}, w_{q+1})$ replaced by $(v^s_{q+1}, w^s_{q+1})$, $R_{\mathrm{corrector}}(v^s_{q+1}, w^s_{q+1})$ the corrector part with $(v_{q+1}, w_{q+1})$ replaced by $(v^s_{q+1}, w^s_{q+1})$, and $R_{\mathrm{oscillation}}(v^s_{q+1}, w^s_{q+1}, R^s_q)$ the oscillation part with $(v_{q+1}, w_{q+1}, R_q)$ replaced by $(v^s_{q+1}, w^s_{q+1}, R^s_q)$.

On the other hand, as indicated in Remark \ref{rk-comparison2}, the increments $(v^s_{q+1}, w^s_{q+1})$ satisfy the same estimates of $(v_{q+1}, w_{q+1})$ in Lemma \ref{le-vq} and Lemma \ref{le-w}. Therefore, $R^s_{\mathrm{linear}}$, $R^s_{\mathrm{corrector}}$ and $R^s_{\mathrm{oscillation}}$ satisfy the same estimates as $R_{\mathrm{linear}}$, $R_{\mathrm{corrector}}$ and $R_{\mathrm{oscillation}}$, respectively. Hence, in view of Lemma \ref{le-R}, we have
\begin{equation}\notag
\|\tilde R^s_{q+1}\|_{L^p}\lesssim \lambda_{q+1}^{-2\varepsilon_R}\delta_{q+2}
\end{equation}
for $p>1$ sufficiently close to $1$ and $\varepsilon_R>0$ sufficiently small.

\cbdu

\begin{Lemma} \label{le-est-mq1e}
For $p>1$ sufficiently close to $1$ and $\varepsilon_R>0$ sufficiently small, the vector $M_{q+1}^\epsilon$ satisfies 
\begin{equation}\notag
\|M_{q+1}^\epsilon\|_{L^p}\leq \lambda_{q+1}^{-\varepsilon_R}\delta_{q+2}+C\|z_{q+1}\|_{L^2},
\end{equation}
with some absolute constant $C>0$.
\end{Lemma}
\pf
An obvious rearrangement yields 
\begin{equation}\notag
\begin{split}
M_{q+1}^\epsilon=B_{q+1}\times u_{q+1}-B_{q}\times u_{q}
= v^s_{q+1}\times u_{q+1}+B_q\times z^s_{q+1}.
\end{split}
\end{equation}
Therefore, we deduce from (\ref{v-pp}), (\ref{v-cp}), (\ref{v-tp}) and (\ref{pp1})-(\ref{pp2}), for $p>1$ close enough to 1,
\begin{equation}\notag
\begin{split}
\|M_{q+1}^\epsilon\|_{L^p}\leq &\ \| v^s_{q+1}\times u_{q+1}\|_{L^p}+\|B_q\times z^s_{q+1}\|_{L^p}\\
\lesssim &\ \| v^s_{q+1}\|_{L^{\frac{2p}{2-p}}} \|u_{q+1}\|_{L^2}+\|B_q\|_{L^{\frac{2p}{2-p}}}\|z^s_{q+1}\|_{L^2}\\
\lesssim &\ \| v^s_{q+1}\|_{L^{\frac{2p}{2-p}}} \|u_{q+1}\|_{L^2}+\|z^s_{q+1}\|_{L^2}\sum_{j=0}^{j=q-1}\|v^s_{j+1}\|_{L^{\frac{2p}{2-p}}}\\
\lesssim &\ \lambda_{q+1}^{-1}\delta_{q+1}^{1/2}r^{\frac32-\frac{3(2-p)}{2p}}
+\|z^s_{q+1}\|_{L^2}\sum_{j=0}^{j=q-1}\lambda_{j+1}^{-1}\delta_{j+1}^{1/2}\lambda_{j+1}^{\frac34\left[\frac32-\frac{3(2-p)}{2p}\right]}\\
\lesssim &\ \lambda_{q+1}^{-1}\delta_{q+1}^{1/2}r^{3-\frac{3}{p}}+\|z^s_{q+1}\|_{L^2}\\
\lesssim &\ \lambda_{q+1}^{-2\varepsilon_R}\delta_{q+2}+\|z^s_{q+1}\|_{L^2}.
\end{split}
\end{equation}

\cbdu

According to Lemma \ref{le-hmhd-new-stress}, Lemma \ref{le-est-tilde-rq1s} and Lemma \ref{le-est-mq1e}, (\ref{RM-h}), (\ref{est-rqs-tilde}) and (\ref{est-increm-M}) in Proposition \ref{hmhd-iterative} are proved to hold for $q+1$.

\medskip

\subsection{Other estimates of Proposition \ref{hmhd-iterative}}\label{sec-other-estimate}
Regarding the energy iteration properties (\ref{energy-q-h}) and (\ref{energy-q1-h}), they can be obtained in a similar way as of (\ref{energy-q}) and (\ref{energy-q1}). Indeed, we notice that the $J_q$ equation in (\ref{hmhd-q}) differs from the $J_q$ equation (\ref{eq-div-q}) by the nonlinear term $\nabla\times\nabla\times (u_q\times B_q)$, which is smaller than the nonlinear portion of the Hall term $\nabla\times\nabla\times (J_q\times B_q)$ (up to scale $\lambda_{q}^{-1}$). Therefore, when deriving (\ref{energy-q-h}) and (\ref{energy-q1-h}) for $q+1$,  the nonlinear term $\nabla\times\nabla\times (u_{q+1}\times B_{q+1})$ can be treated as a small error term and hence absorbed by other terms in the estimates. Thus, slight modification of the proof of energy iteration in \cite{BV} will yield (\ref{energy-q-h}) and (\ref{energy-q1-h}).

We conclude the proof of Proposition \ref{hmhd-iterative}. 

\medskip

\subsection{Proof of Theorem \ref{thm}}
We are left to show that the sequence $\{(u_q, B_q)\}_{q=1}^{\infty}$ converges to a pair $(u,B)\in \left(L^\infty(0,T; L^2(\mathbb T^3))\cap L^2(0,T; H^1(\mathbb T^3))\right)^2$ which solves the Hall MHD (\ref{HMHD}).

For given $(u_0, B_0, J_0, R_0)$, we apply Proposition \ref{hmhd-iterative} iteratively to obtain a sequence of approximating solutions $\{(u_q, B_q, J_q, R_q)\}$ satisfying (\ref{bq-c1-h})-(\ref{est-increm-M}). It follows from
(\ref{est-increm-h}) that  
\begin{equation}\notag
\sum_{q\geq0}\|J_{q+1}-J_q\|_{L^2}= \sum_{q\geq0}\|w^s_{q+1}\|_{L^2}\lesssim \sum_{q\geq0}\delta_{q+1}^{1/2}<\infty.
\end{equation}
which implies the strong convergence of $J_q=\nabla\times B_q$ to a function $J$ in $C^0(0,T;L^2)$, and the strong convergence of $B_q$ to a function $B$ in $C^0(0,T;H^1)$ with $J=\nabla\times B$ and $\nabla\cdot B=0$. 

According to the analysis above, 
we have $u_{q}\in L^\infty(0,T; L^2(\mathbb T^3))\cap L^2(0,T; H^1(\mathbb T^3))$ and $B_q\in C^0(0,T;H^1)$ for all $q\geq 1$. We claim that,  there exists a subsequence (with the same notation) such that
\begin{equation}\label{convergence-uq}
u_q\to u \ \ \mbox{weakly in} \ L^2(0,T; H^1(\mathbb T^3)) \ \mbox{and strongly in} \ L^2(0,T; L^2(\mathbb T^3)).
\end{equation}
The weak convergence in $L^2(0,T; H^1(\mathbb T^3))$ is automatic. In order to show the strong convergence in $L^2(0,T; L^2(\mathbb T^3))$, we apply the Aubin-Lions lemma (cf. \cite{Tem77}). Since $u_q\in L^2(0,T; H^1(\mathbb T^3))$, $H^1$ is compactly embedded in $L^2$ and $L^2$ is continuously embedded in $H^{-1}$, it is sufficient to prove 
\begin{equation}\label{H-negative}
\partial_t u_q\in L^{\frac43}([0,T]; H^{-1}(\mathbb T^3)). 
\end{equation}
Hence, the Aubin-Lions lemma implies the strong convergence of $u_q$ in $L^2(0,T; L^2(\mathbb T^3))$. We prove (\ref{H-negative}) as follows. Let $\varphi\in H^1(\mathbb T^3)$ with $\nabla\cdot \varphi=0$. Taking inner product of equation (\ref{nse-q}) with $\varphi$ (to simplify notations, $q+1$ in (\ref{nse-q}) is replaced by $q$) and applying integration by parts yields
\begin{equation}\notag
\int_{\mathbb T^3} \partial_t u_q\cdot \varphi\, dx= \int_{\mathbb T^3} (u_q\cdot\nabla) \varphi \cdot u_q\, dx-\int_{\mathbb T^3}\nabla u_q : \nabla \varphi\, dx-\int_{\mathbb T^3} (B_q\cdot\nabla) \varphi \cdot B_q\, dx.
\end{equation}
Thus, we obtain
\begin{equation}\notag
\left|\int_{\mathbb T^3} \partial_t u_q\cdot \varphi\, dx\right|\leq C \|\varphi\|_{H^1(\mathbb T^3)}\left(\|u_q\|_{L^4(\mathbb T^3)}^2+\|B_q\|_{L^4(\mathbb T^3)}^2+\|\nabla u_q\|_{L^2(\mathbb T^3)}\right)
\end{equation}
for a constant $C>0$, 
which implies
\begin{equation}\label{H-negative1}
\|\partial_t u_q\|_{H^{-1}(\mathbb T^3)}\leq C \left(\|u_q\|_{L^4(\mathbb T^3)}^2+\|B_q\|_{L^4(\mathbb T^3)}^2+\|\nabla u_q\|_{L^2(\mathbb T^3)}\right).
\end{equation}
In view of Gagliardo-Nirenberg's inequality, we have 
\begin{equation}\label{H-negative2}
\begin{split}
\|u_q\|_{L^4(\mathbb T^3)}^2\leq&\ 4 \|u_q\|_{L^2(\mathbb T^3)}^{\frac12} \|\nabla u_q\|_{L^2(\mathbb T^3)}^{\frac32},\\
\|B_q\|_{L^4(\mathbb T^3)}^2\leq&\ 4 \|B_q\|_{L^2(\mathbb T^3)}^{\frac12} \|\nabla B_q\|_{L^2(\mathbb T^3)}^{\frac32}.
\end{split}
\end{equation} 
Applying the fact of $u_q\in L^2(0,T; L^2(\mathbb T^3))$ and $B_q\in C^0(0,T;H^1)$, and putting (\ref{H-negative1}) and (\ref{H-negative2}) together gives
\begin{equation}\notag
\|\partial_t u_q\|_{H^{-1}(\mathbb T^3)}\leq C \left(\|u_q\|_{H^1(\mathbb T^3)}^{\frac32}+\|u_q\|_{H^1(\mathbb T^3)}+1\right)
\end{equation}
for a constant $C>0$ that is independent of $q$. As a consequence, we deduce
\begin{equation}\notag
\begin{split}
&\int_0^T\|\partial_t u_q\|_{H^{-1}(\mathbb T^3)}^{\frac43} \, dt\\
\leq&\  C \int_0^T\|u_q\|_{H^{1}(\mathbb T^3)}^{\frac43} \, dt+C \int_0^T\|u_q\|_{H^{1}(\mathbb T^3)}^2 \, dt+CT\\
\leq &\ C T^{\frac13} \left(\int_0^T\|u_q\|_{H^{1}(\mathbb T^3)}^2 \, dt\right)^{\frac23}
+C \int_0^T\|u_q\|_{H^{1}(\mathbb T^3)}^2 \, dt+CT\\
\leq &\ C(T)
\end{split}
\end{equation}
for a constant $C(T)$ dependent on $T$ and independent on $q$, where we used the fact $u_q\in L^2(0,T; H^1(\mathbb T^3))$. This completes the proof of (\ref{H-negative}) and hence (\ref{convergence-uq}).

Combining (\ref{convergence-uq}) and the fact of $B_q$ converging to $B$ strongly in $C^0(0,T;H^1(\mathbb T^3))$, we conclude that $(u,B)$ solves the NSE part of (\ref{HMHD}) in the weak sense.

On the other hand, the facts $\|R_q\|_{L^\infty(0,T; L^1(\mathbb T^3))}\to 0$ and $\|M_q^{\epsilon}\|_{L^\infty(0,T; L^1(\mathbb T^3))}\to 0$ as $q\to \infty$ lead to $\|R_q^s\|_{L^\infty(0,T; L^1(\mathbb T^3))}\to 0$ as $q\to \infty$. Thus, $(u,B)$ also solves the second equation of (\ref{HMHD}) in the weak sense. It follows that $(u, B)$ is a weak solution of system  (\ref{HMHD}).

\bigskip

\section{Appendix: Vector calculus identities}
\label{sec-vec}

Let $A$ and $B$ be vector valued functions. The following identities hold:
\begin{equation}\label{identities-vec}
\begin{split}
\nabla\times(A\times B)=&\ A(\nabla\cdot B)-B(\nabla\cdot A)+(B\cdot\nabla) A-(A\cdot \nabla)B;\\
(A\cdot\nabla)B=&\ A\nabla B-A\times (\nabla\times B), \ \ \mbox{with} \ \ A\nabla B=(A_i\partial_j B_i)_{j};\\
\nabla\times(\nabla\times A)=&\ \nabla(\nabla\cdot A)-\nabla^2A=\nabla(\nabla\cdot A)-\Delta A.
\end{split}
\end{equation}
Assume $B$ is divergence free. Let $J=\nabla\times B$, then $\nabla\cdot J=0$ and $\nabla\cdot (\nabla\times J)=0$.
Applying the identities of (\ref{identities-vec}), one can rewrite 
\begin{equation}\notag
\begin{split}
&\nabla\times\nabla\times[(\nabla\times B)\times B]\\
=&\ \nabla\times((B\cdot\nabla) J-(J\cdot\nabla)B)\\
=&\ \nabla\times (B\nabla J)-\nabla\times (B\times (\nabla\times J))-\nabla\times (J\nabla B)+\nabla\times (J\times (\nabla\times B))\\
=& -2\nabla\times (J\nabla B)+(B\cdot\nabla)(\nabla\times J)-((\nabla\times J)\cdot\nabla)B.
\end{split}
\end{equation}
One can also derive the identity,
\[\Delta(\nabla\times B)=\nabla\times (\Delta B).\]
Thus, taking curl of the Hall equation
\begin{equation}\notag
B_t+\nabla\times((\nabla\times B)\times B)=\Delta B,
\end{equation} 
we obtain the equation of the current density $J=\nabla\times B$,
\begin{equation}\notag
J_t+\nabla\cdot[B\otimes (\nabla\times J)-(\nabla\times J)\otimes B] -2\nabla\times (J\nabla B)=\Delta J.
\end{equation}

In general, for $A$ and $B$ with $\nabla\cdot A=\nabla\cdot B=0$, we have
\begin{equation}\notag
\begin{split}
&\nabla\times\nabla\times (A\times B)\\
=&\ \nabla\cdot[B\otimes (\nabla\times A)-(\nabla\times A)\otimes B]+\nabla\cdot[(\nabla\times B)\otimes A-A\otimes (\nabla\times B)]\\
& -2\nabla\times (A\nabla B).
\end{split}
\end{equation}

\bigskip

\section*{Acknowledgement}

The author is indebted to the anonymous referees for their valuable suggestions and comments, which have helped improve the manuscript a lot.


\bigskip


\begin{thebibliography}{XX}


\bibitem{ADFL}
M. Acheritogaray, P. Degond, A. Frouvelle and J-G. Liu.
\newblock {\em Kinetic formulation and global existence for the Hall-Magnetohydrodynamic system}.
\newblock Kinetic and Related Models, 4: 901--918, 2011.

\bibitem{BBV}
R. Beekie, T. Buckmaster, and V. Vicol.
\newblock {\em Weak solutions of ideal MHD which do not conserve magnetic helicity}.
\newblock Ann. of PDE, Vol. 6, No. 1, 2020.

\bibitem{BCV}
T. Buckmaster, M. Colombo, and V. Vicol.
\newblock {\em Wild solutions of the Navier-Stokes equations whose singular sets in time have Hausdorff dimension strictly less than 1}.
\newblock Journal of the EMS, in press.

\bibitem{BDLIS}
T. Buckmaster, C. De Lellis, P. Isett, and L. Sz\'ekelyhidi.
\newblock {\em Anomalous dissipation for $1/5$-H\"older Euler flows}.
\newblock Ann. of Math., Vol. 182, No. 1: 127-172, 2015.

\bibitem{BDLS}
T. Buckmaster, C. De Lellis, and L. Sz\'ekelyhidi.
\newblock {\em Dissipative Euler flows with Onsager-critical spatial regularity}.
\newblock Comm. Pure Appl. Math., Vol. 69 No. 9, 16131670, 2016.

\bibitem{BDLSV}
T. Buckmaster, C. De Lellis, L. Sz\'ekelyhidi, and V. Vicol.
\newblock {\em Onsager's conjecture for admissible weak solutions}.
\newblock Comm. Pure Appl. Math., https://doi.org/10.1002/cpa.21781. 2018.

\bibitem{BSV}
T. Buckmaster, S. Shkoller, and V. Vicol.
\newblock {\em Nonuniqueness of weak solutions to the SQG euqation}.
\newblock Comm. Pure Appl. Math., 2019(72):1809--1874, 2019.

\bibitem{BV19}
T. Buckmaster, and V. Vicol.
\newblock {\em Convex integration and phenomenologies in turbulence}.
\newblock EMS Surveys in Mathematical Sciences, 6(1):173--263, 2020.

\bibitem{BV21}
T. Buckmaster, and V. Vicol.
\newblock {\em Convex integration constructions in hydrodynamics}.
\newblock Bulletin of the AMS, 58(1): 1--44, 2021.

\bibitem{BV}
T. Buckmaster, and V. Vicol.
\newblock {\em Nonuniqueness of weak solutions to the Navier-Stokes equation}.
\newblock Ann. of Math., Vol. 189, No. 1: 101-144, 2019.

\bibitem{BMS}
J. Burczak, S. Modena, and L. Sz\'ekelyhidi .
\newblock {\em Non uniqueness of power-law flows}.
\newblock arXiv: 2007.08011, 2020.

\bibitem{CDL}
D. Chae, P. Degond and J-G. Liu.
\newblock {\em Well-posedness for Hall-magnetohydrodynamics}.
\newblock Ann. Inst. H. Poincar\'e Anal. Non Lineaire, Vol. 31: 555--565, 2014.

\bibitem{CL}
D. Chae and J. Lee.
\newblock {\em On the blow-up criterion and small data global existence for the Hall-magneto-hydrodynamics}.
\newblock J. Differential Equations, 256: 3835--3858, 2014.

\bibitem{CS}
D. Chae,  and M. Schonbek.
\newblock {\em On the temporal decay for the Hall-magnetohydrodynamic equations}.
\newblock J. Differential Equations,  Vol. 255: 3971--3982, 2013.

\bibitem{CWW}
D. Chae,  R. Wan and J. Wu.
\newblock {\em Local well-posedness for the Hall--MHD equations with fractional magnetic diffusion}.
\newblock arXiv:1404.0486v2, 2014.

\bibitem{CWeng}
D. Chae and S. Weng.
\newblock {\em Singularity formation for the incompressible Hall-MHD equations without resistivity}.
\newblock Ann. I. H. Poincar\'e-AN, Vol. 33: 1009--1022, 2016.

\bibitem{CW}
D. Chae and J. Wolf.
\newblock {\em On partial regularity for the 3D non-stationary Hall magnetohydrodynamics equations on the plane}.
\newblock Comm. Math. Phys., Vol. 354: 213--230, 2017.

\bibitem{CL1}
A. Cheskidov and X. Luo.
\newblock {\em Anomalous dissipation, anomalous work, and energy balance for smooth solutions of the Navier-Stokes equations}.
\newblock arXiv: 1910.04204, 2019.

\bibitem{CL2}
A. Cheskidov and X. Luo.
\newblock {\em $L^2$-critical nonuniqueness for the 2D Navier-Stokes equations}.
\newblock arXiv: 2105.12117, 2021.

\bibitem{CL3}
A. Cheskidov and X. Luo.
\newblock {\em Sharp nonuniqueness for the Navier-Stokes equations}.
\newblock arXiv:2009.06596, 2020.


\bibitem{CDLDR}
M. Colombo, C. De Lellis, and L. De Rosa.
\newblock {\em Ill-posedness of Leray solutions for the hypodissipative Navier-Stokes equations}.
\newblock Comm. Math. Phys., Vol.362, No.2, 659688, 2018.



\bibitem{Dai1}
M. Dai.
\newblock {\em Local well-posedness of the Hall-MHD system in $H^s(\mathbb R^n)$ with $s>\frac n2$}.
\newblock Mathemaische Nachrichten, Vol. 293, Iss. 1: 67--78.

\bibitem{Dai2}
M. Dai.
\newblock {\em Local well-posedness for the Hall-MHD system in optimal Sobolev spaces}.
\newblock arXiv: 1803.09556, 2018.



\bibitem{Dai4}
M. Dai.
\newblock {\em Regularity criterion for the 3D Hall-magneto-hydrodynamics}.
\newblock Journal of Differential Equations. Vol. 261: 573--591, 2016.



\bibitem{DL}
M. Dai and H. Liu.
\newblock {\em Long time behavior of solutions to the 3D Hall-magneto-hydrodynamics system with one diffusion}.
\newblock Journal of Differential Equations, Vol. 266: 7658--7677, 2019.

\bibitem{DSz}
S. Daneri, and L. Sz\'ekelyhidi.
\newblock {\em Non-uniqueness and h-principle for H\"older-continuous weak solutions of the Euler equations}.
\newblock Arch. Ration. Mech. Anal., Vol. 224 No.2: 471--514, 2017.

\bibitem{DLS1}
C. De Lellis, and L. Sz\'ekelyhidi.
\newblock {\em Dissipative continuous Euler flows}.
\newblock Invent. Math., Vol.193 No. 2: 377--407, 2013.

\bibitem{DLS3}
C. De Lellis, and L. Sz\'ekelyhidi, Jr.
\newblock {\em Dissipative Euler flows and Onsager's conjecture}.
\newblock Journal of the European Mathematical Society,  16(7): 1467--1505, 2014.

\bibitem{DLS2}
C. De Lellis, and L. Sz\'ekelyhidi.
\newblock {\em The Euler equations as a differential inclusion}.
\newblock Ann. of Math.,  Vol.170 No.3: 1417--1436, 2009.

\bibitem{DR}
L. De Rosa.
\newblock {\em Infinitely many Leray–Hopf solutions for the fractional Navier-Stokes equations}.
\newblock Communications in Partial Differential Equations, 44(4): 335--365, 2019.

\bibitem{DS}
E. Dumas and F. Sueur.
\newblock {\em On the weak solutions to the Maxwell-Landau-Lifshitz equations and to the Hall-magnetohydrodynamic equations}.
\newblock Comm. Math. Phys., 330: 1179--1225, 2014.

\bibitem{FLS}
D. Faraco, S. Lindberg, and L. Sz\'ekelyhidi.
\newblock {\em Bounded solutions of ideal MHD with compact support in space-time}.
\newblock Arch. Ration. Mech. Anal., https://doi.org/10.1007/s00205-020-01570-y, 2020.

\bibitem{FMRR}
C. L. Fefferman, D. S. McCormick, J. C. Robinson and J. L. Rodrigo.
\newblock {\em Higher order commutator estimates and local existence for the non-resistive MHD equations and related models}.
\newblock Journal of Functional Analysis, Vol. 267: 1035--1056, 2014.


\bibitem{Is0}
P. Isett.
\newblock {\em Holder continuous Euler flows with compact support in time}.
\newblock Ph. D. Thesis, Princeton Jniversity, 2013.

\bibitem{Is}
P. Isett.
\newblock {\em A Proof of Onsager's Conjecture}.
\newblock Ann. of Math., Vol.188 No.3: 1--93, 2018.

\bibitem{JS1}
H. Jia and V. \v{S}ver\'ak.
\newblock {\em Are the incompressible 3d Navier-Stokes equations locally ill-posed in the natural energy space?}
\newblock J. Funct. Anal., Vol. 268(12): 3734--3766, 2015.

\bibitem{JS2}
H. Jia and V. \v{S}ver\'ak.
\newblock {\em Local-in-space estimates near initial time for weak solutions of the Navier-Stokes equations and forward self-similar solutions}.
\newblock Invent. Math., Vol. 196(1): 233--265, 2014.



\bibitem{Lions}
J.-L. Lions.
\newblock {\em Quelques m\'ethodes de r\'esolution des probl\'emes aux limites non lin\'eaires}.
\newblock volume 1. Dunod; Gauthier-Villars, Paris, 1969.



\bibitem{LT}
T. Luo and E. Titi.
\newblock {\em Non-uniqueness of weak solutions to hyperviscous Navier-Stokes equations-On sharpness of J.-L. Lions exponent}.
\newblock arXiv:1808.07595, 2018.


\bibitem{MS}
S. Modena, and L. Sz\'ekelyhidi.
\newblock {\em Non-uniqueness for the transport equation with Sobolev vector fields}.
\newblock Annals of PDE, Vol. 4 No. 18,  2018. https://doi.org/10.1007/s40818-018-0056-x

\bibitem{MSv}
S. M\"uller, and V. \v{S}ver\'ak..
\newblock {\em Convex integration for Lipschitz mappings and counterexamples to regularity}.
\newblock Ann. of Math., 157(3): 715--742, 2003.

\bibitem{No}
M. Novack.
\newblock {\em Nonuniqueness of weak solutions to the 3 dimensional quasi-geostrophic equations}.
\newblock SIAM J. Math. Anal., 52(4): 3301--3349, 2020.


\bibitem{On}
L. Onsager.
\newblock {\em Statistical hydrodynamics}.
\newblock Nuovo Cimento, no. Suppl, Supplemento, no. 2 (Convegno Internazionale di Meccanica Statistica), 9(6): 279--287, 1949.

\bibitem{Tao}
T. Tao.
\newblock {\em Finite time blowup for an averaged three-dimensional Navier-Stokes equation}.
\newblock J. Amer. Math. Soc., Vol. 29: 601--674, 2016.

\bibitem{Tem77}
R. Temam.
\newblock {\em Navier-Stokes Equations, rev. ed., Studies in Mathematics and its Applications 2}.
\newblock North-Holland, Amsterdam, 1977.

\bibitem{Tem}
R. Temam.
\newblock {\em Navier-Stokes Equations: Theory and Numerical Analysis}.
\newblock AMS Chelsea, Providence, Rhode Island, 2000.

\bibitem{Sch}
V. Scheffer.
\newblock {\em An inviscid flow with compact support in space-time}.
\newblock J. Geom. Anal., 3(4):343--401, 1993.


\bibitem{Sh1}
A. Shnirelman.
\newblock {\em On the nonuniqueness of weak solution of the Euler equations}.
\newblock Comm. Pure Appl. Math., 50(12):1261--1286, 1997.

\bibitem{Sh2}
A. Shnirelman.
\newblock {\em Weak solutions with decreasing energy of incompressible Euler equations}.
\newblock Comm. Math. Phys., 210(3):541--603, 2000.


\end{thebibliography}
\end{document}